\documentclass[preprint,11pt,numbers]{elsarticle-arxiv}
\usepackage{lmodern}
\usepackage{amssymb,amsmath,amsthm,mathtools}
\usepackage{booktabs,multirow,enumitem,indentfirst,etoolbox,float,url}
\makeatletter
\IfFileExists{placeins.sty}{\usepackage{placeins}}{%
  \newcommand{\FloatBarrier}{\par\ifx\@deferlist\@empty\else\clearpage\fi}}
\makeatother
\usepackage{tikz}
\usetikzlibrary{arrows.meta}
\tikzset{
  graph vertex/.style={circle,fill=black,inner sep=1.7pt},
  marked edge/.style={draw=blue!65!black,line width=1.05pt},
  unmarked edge/.style={draw=black!65,line width=0.65pt}
}
\usepackage[hidelinks]{hyperref}
\ifdefined\pdfgentounicode
  \InputIfFileExists{glyphtounicode.tex}{}{}
\fi
\setcitestyle{square,comma}
\numberwithin{equation}{section}
\newtheorem{theorem}{Theorem}[section]
\newtheorem{conjecture}[theorem]{Conjecture}
\newtheorem{lemma}[theorem]{Lemma}
\theoremstyle{definition}

\newenvironment{lateproof}[1]{\begin{proof}[Proof of #1]}{\end{proof}}
\newcommand{\papertitle}{The Proof of the Faudree--Gy\'arf\'as--Schelp--Tuza Conjecture}
\hypersetup{pdftitle={\papertitle},pdfauthor={Zhijun Lu, Chengshuo Qiao, Qirui Ying, Huimin Song}}
\setlist[enumerate]{itemsep=3pt,topsep=5pt,parsep=0pt,listparindent=1.5em}
\patchcmd{\abstract}
  {\par\medskip\noindent\unskip\ignorespaces}
  {\par\medskip\ignorespaces}
  {}{\PackageError{manuscript}{Could not apply abstract indentation}{}}

\begin{document}
\bibliographystyle{alpha}
\begin{frontmatter}
\title{\textbf{\Large\papertitle}}
\author[1]{Zhijun Lu\fnref{equalcontrib}}
\author[1]{Chengshuo Qiao\fnref{equalcontrib}}
\author[1]{Qirui Ying\fnref{equalcontrib}}
\author[1]{Huimin Song\corref{cor1}}
\fntext[equalcontrib]{Zhijun Lu, Chengshuo Qiao and Qirui Ying contributed equally to this work.}
\cortext[cor1]{Corresponding author. \textit{E-mail address:}
\url{hmsong@sdu.edu.cn} (H.Song).}
\affiliation[1]{organization={School of Mathematics and Statistics},
 addressline={Shandong University},city={Weihai},postcode={264209},
 state={Shandong},country={China}}

\begin{abstract}
In 1990, Faudree, Gy\'arf\'as, Schelp and Tuza proposed a relaxation
of their multigraph version of the Erd\H{o}s--Ne\v{s}et\v{r}il conjecture
on strong edge coloring. They conjectured that every
finite loopless multigraph $G$ satisfies
$\omega(L(G)^2)\leq\Delta^2+\lfloor\frac{\Delta}{2}\rfloor^2$, where
$\Delta=\Delta(G)$ is the maximum degree, $L(G)$ is the line graph
of $G$, and $\omega$ is the clique number. In this paper, we prove
this conjecture. We first introduce the weighted flag algebra method,
which generalizes the local flag algebra method of
Davey et al.~\href{https://arxiv.org/abs/2607.12461}{[arXiv:2607.12461]},
to show that $\omega(L(G)^2)\leq\frac{5}{4}\Delta^2$ for every finite loopless
multigraph $G$. This settles the even-degree case.
For odd $\Delta\geq 9$, we show that if
$\omega(L(G)^2)>\Delta^2+\lfloor\frac{\Delta}{2}\rfloor^2$,
then $L(G)^2$ has a clique exceeding this bound whose vertices
correspond to a set of edges of $G$ that includes both a pair of
parallel edges and five edges forming an induced five-cycle in the
underlying simple graph of $G$. Counting around such cycles gives a
bound that contradicts the weighted flag inequalities. The remaining
small degrees are handled separately, thereby completing the proof.
\end{abstract}
\end{frontmatter}

\thispagestyle{plain}
\tableofcontents
\clearpage
\section{Introduction}\label{sec:intro}

\subsection{Preliminaries}\label{sec:prelim}

All graphs are finite and simple unless stated otherwise. In a loopless
multigraph, parallel edges are distinct and contribute to degrees with
multiplicity. For a graph $G$, write $V(G)$, $E(G)$ and $\Delta(G)$ for
its vertex set, edge set and maximum degree. For $v\in V(G)$, let
$N_G(v)$ be its set of neighbors and let $d_G(v)$ be its degree.
For $X\subseteq V(G)$, write $G[X]$ for the subgraph induced by $X$ and
$d_G(v,X)$ for the number of edges from $v$ to $X$.
A graph is $d$-regular if all its vertices have degree $d$. Write $C_5$
for the five-cycle and $K_n$ for the complete graph on $n$ vertices.

A vertex cover of $G$ is a set meeting every edge of $G$. An independent
set is a set of pairwise nonadjacent vertices. For a subgraph $H$ of $G$
with at least one edge, its Ore-degree in $G$ is
$\max_{uv\in E(H)}(d_G(u)+d_G(v))$.

The \emph{underlying simple graph} of a loopless multigraph $G$
has vertex set $V(G)$ and one edge
$uv$ for each pair of distinct vertices $u,v$ joined by an edge of $G$.
The line graph $L(G)$ has
one vertex for each edge of $G$, with two distinct edges adjacent
precisely when they share an endpoint in $G$. For a graph $H$,
its square $H^2$ has vertex set $V(H)$, with two distinct vertices
adjacent precisely when $H$ contains a path of length at most two
between them. We write $\omega(H)$ for the maximum size of a clique in
$H$, that is, a set of pairwise adjacent vertices. A set
$S\subseteq E(G)$ is a \emph{strong clique} in $G$ if it is a clique
in $L(G)^2$.

For a positive integer $k$, write $[k]=\{1,\ldots,k\}$. A proper vertex
$k$-coloring of $H$ is a map $V(H)\to[k]$ assigning distinct colors to
adjacent vertices. The chromatic number $\chi(H)$ is the least $k$ for
which $H$ admits such a coloring, with $\chi(H)=0$ when $V(H)$ is empty.
A proper vertex coloring of $L(G)^2$ gives a \emph{strong edge
coloring} of $G$. The minimum number of colors in such a coloring is
the \emph{strong chromatic index} $\chi(L(G)^2)$.

A real symmetric matrix $M$ is positive semidefinite, written
$M\succeq0$, if $x^{\mathsf T}Mx\geq0$ for every real column vector
$x$ of the appropriate size.

\subsection{Background and Main Result}\label{sec:background}

In 1985, at a seminar in Prague, Erd\H{o}s and Ne\v{s}et\v{r}il asked
how large the strong chromatic index can be for a graph of given
maximum degree~\cite{FGST1989,ErdosNesetril1989}. Erd\H{o}s~\cite{Erdos1988}
stated their conjecture with the bound $\frac{5}{4}\Delta(G)^2$.
Faudree, Schelp, Gy\'arf\'as and Tuza~\cite{FGST1990} formulated a
version for loopless multigraphs that gives a sharper bound for odd
maximum degree. For a nonnegative integer $d$, let
$f(d)=d^2+\lfloor\frac{d}{2}\rfloor^2$. Thus $f(d)=\frac{5d^2}{4}$
when $d$ is even and $f(d)=\frac{5d^2-2d+1}{4}$ when $d$ is odd.

\begin{conjecture}[\cite{FGST1990}]
\label{conj:strong-edge-coloring}
Every finite loopless multigraph $G$ satisfies
$\chi(L(G)^2)\leq f(\Delta(G))$.
\end{conjecture}

Faudree, Gy\'arf\'as, Schelp and Tuza then asked for the largest size of
a set of edges that must receive pairwise distinct colors in every
strong edge coloring. They proposed the following relaxation of
\textbf{Conjecture~\ref{conj:strong-edge-coloring}}.

\begin{conjecture}[\cite{FGST1990}]\label{conj:fgst}
Every finite loopless multigraph $G$ satisfies
$\omega(L(G)^2)\leq f(\Delta(G))$.
\end{conjecture}

\textbf{Table~\ref{tab:strong-clique-progress}} summarizes the progress
on \textbf{Conjecture~\ref{conj:fgst}} for simple graphs and loopless
multigraphs. Here $\Delta=\Delta(G)$ and $\varepsilon>0$ is an absolute
constant.

\begin{table}[H]
\centering
\makeatletter
\renewcommand{\@makecaption}[2]{\centering #1: #2\par\smallskip}
\makeatother
\caption{Previous upper bounds on $\omega(L(G)^2)$.}
\label{tab:strong-clique-progress}
\renewcommand{\arraystretch}{1.2}
\setlength{\tabcolsep}{5pt}
\setlength{\arrayrulewidth}{0.35pt}
\begin{tabular}{|c|c|c|c|}
\hline
Graph class & Authors & Upper bound & Reference\\
\hline
{\renewcommand{\arraystretch}{1}%
  \begin{tabular}[c]{@{}c@{}}Loopless\\multigraphs\end{tabular}}
  & {\renewcommand{\arraystretch}{1}%
     \begin{tabular}[c]{@{}c@{}}Faudree, Schelp,\\Gy\'arf\'as and Tuza\end{tabular}}
  & $\omega(L(G)^2)\leq(2-\varepsilon)\Delta^2$
  & \cite{FGST1990}\\
\hline
\multirow{4}{*}{Simple graphs}
  & Bruhn and Joos
  & $\omega(L(G)^2)\leq1.74\Delta^2$ ($\Delta\geq400$)
  & \cite{BruhnJoos2018}\\
\cline{2-4}
  & \'{S}leszy\'{n}ska-Nowak
  & $\omega(L(G)^2)\leq\frac{3}{2}\Delta^2$
  & \cite{SleszynskaNowak2016}\\
\cline{2-4}
  & Faron and Postle
  & $\omega(L(G)^2)\leq\frac{4}{3}\Delta^2$
  & \cite{FaronPostle2019}\\
\cline{2-4}
  & Kumar, Mohar and Pragada
  & $\omega(L(G)^2)\leq\frac{2607}{1987}\Delta^2$
  & \cite{KumarMoharPragada2026}\\
\hline
\end{tabular}
\end{table}

In this paper, we prove \textbf{Conjecture~\ref{conj:fgst}}.

\begin{theorem}\label{thm:multigraph}
Every finite loopless multigraph $G$ satisfies
$\omega(L(G)^2)\leq f(\Delta(G))$.
\end{theorem}

\begin{figure}[!htbp]
\centering
\begin{tikzpicture}[scale=0.9,every node/.style={font=\small}]
\foreach \i in {0,...,4}
  \coordinate (v\i) at ({90+72*\i}:1.16);
\foreach \i/\j in {0/1,1/2,2/3,3/4,4/0}
  \draw[unmarked edge] (v\i)--(v\j);
\foreach \i in {0,...,4}
  \node[graph vertex] at (v\i) {};
\node at (0,-2.08) {$C_5$};

\begin{scope}[xshift=5.45cm]
\foreach \i in {0,...,4}{
  \pgfmathsetmacro{\angle}{90+72*\i}
  \coordinate (a\i) at (\angle:1.18);
  \coordinate (b\i) at (\angle:1.76);
  \coordinate (part\i) at (\angle:1.47);
  \begin{scope}[shift={(part\i)},rotate=\angle]
    \draw[draw=black!45,fill=black!7] (0,0)
      ellipse[x radius=.42,y radius=.17];
  \end{scope}
}
\foreach \i/\j in {0/1,1/2,2/3,3/4,4/0}{
  \draw[marked edge] (a\i)--(a\j) (a\i)--(b\j)
    (b\i)--(a\j) (b\i)--(b\j);
}
\foreach \i in {0,...,4}{
  \node[graph vertex] at (a\i) {};
  \node[graph vertex] at (b\i) {};
}
\node at (0,-2.08) {$t=2$};
\end{scope}
\end{tikzpicture}
\caption{The cycle $C_5$ and its blow-up with two vertices per part.}
\label{fig:blowup}
\end{figure}
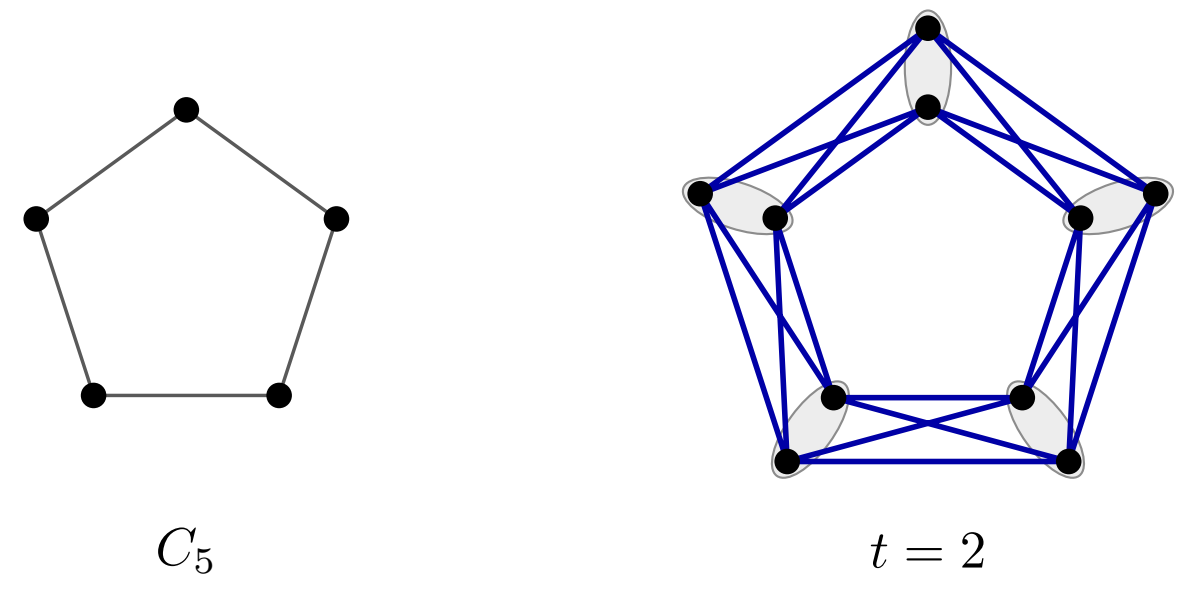

The five-cycle blow-ups of Faudree et al.~\cite{FGST1990}
attain the bound. Let $C_5=v_0v_1v_2v_3v_4v_0$. Replace each $v_i$ by an independent set $V_i$ and add every edge between $V_i$ and $V_{i+1}$, with indices taken modulo five. It is easy to see that the entire edge set of the blow-up is a strong clique. Taking $|V_i|=t$ for each $i$ gives
maximum degree $2t$ and $5t^2=f(2t)$ edges, as illustrated for
$t=2$ in \textbf{Figure~\ref{fig:blowup}}. Taking $|V_2|=|V_3|=r+1$ and
$|V_0|=|V_1|=|V_4|=r$ gives maximum degree $2r+1$ and
$5r^2+4r+1=f(2r+1)$ edges, as required.

\relax\space

\section{The Method of the Proof}\label{sec:overview}

\subsection{The Weighted Flag Algebra Method}\label{subsec:weighted-sketch}

In 2007, Alexander A. Razborov introduced flag algebras to study
asymptotic relations between induced subgraph densities~\cite{Razborov2007}.
For graphs $F$ and $G$ with $|V(G)|\geq |V(F)|$, the
\emph{induced subgraph density} $p(F;G)$ is the probability that
$G[U]$ is isomorphic to $F$, where $U$ is chosen uniformly from
the $|V(F)|$-element subsets of $V(G)$. As a first application,
Razborov gave a new proof of Fisher's lower bound on $p(K_3;G)$
for a prescribed value of $p(K_2;G)$. For an introduction to flag algebras, see
de~Carli~Silva, de~Oliveira~Filho and Sato~\cite{CarliSilvaOliveiraSato2016}.

In July 2026, Davey, Hurley, de Joannis de Verclos, Kang and Volec
introduced local flag algebras~\cite{DaveyEtAlLocal2026}. They used
the method to improve the known upper bound on the strong chromatic
index of simple graphs with sufficiently large maximum degree,
making progress toward
\textbf{Conjecture~\ref{conj:strong-edge-coloring}}~\cite{DaveyEtAlStrong2026}.

In this paper, we use the flag algebra method to bound the size of
a strong clique. We fix a strong clique and mark every edge in it.
In a loopless multigraph, several marked edges may join the same pair
of vertices. If we count each pair joined by a marked edge once, its
contribution is the same regardless of how many marked edges join it,
so the count loses information about their multiplicity. We therefore
introduce the \emph{weighted flag algebra method}, which extends the
local flag algebra method to loopless multigraphs.
To formulate our method, we first introduce some notations and
definitions.

Let $D$ be a positive integer and let $G$ be a finite loopless
multigraph of maximum degree at most $D$. Fix a strong clique of
$m>0$ edges and mark its edges. Delete every vertex that is not an
endpoint of a marked edge, together with its incident edges. For each
adjacent pair $u,v$ in the resulting multigraph, retain all marked
edges with endpoints $u,v$ if there are any, and retain one unmarked
edge otherwise. The reduction preserves every marked edge and every
adjacency between endpoints of marked edges, so the marked edges
still form a strong clique.

For the counts below, we adapt the regularization lemma of
Davey et al.~\cite{DaveyEtAlStrong2026} to embed the reduced
multigraph as an induced subgraph of a $D$-regular loopless
multigraph.

\begin{lemma}\label{lem:regularization}
Let $D\geq1$ be an integer and let $G$ be a finite loopless
multigraph with $\Delta(G)\leq D$. There is a finite $D$-regular
loopless multigraph $\widehat G$ containing $G$ as an induced
subgraph, and every edge added to $G$ has multiplicity one. In
particular, $\widehat G$ is simple if $G$ is simple.
\end{lemma}

\begin{proof}
The sum of the degree deficits in $G$ is
\[
 \sum_{v\in V(G)}\bigl(D-d_G(v)\bigr)=D|V(G)|-2|E(G)|.
\]
If this sum is even, set $G_0=G$. If it is odd, then $D$ is odd,
and we form $G_0$ by adjoining one isolated vertex to $G$. This
increases the sum of the deficits by $D$, so the sum for $G_0$
is even in either case.

Define
\[
 \mathcal I=\bigl\{(v,i):v\in V(G_0),\ i\in\mathbb Z,\
                     1\leq i\leq D-d_{G_0}(v)\bigr\}.
\]
Each $(v,i)\in\mathcal I$ indexes one of the edges to be added
at $v$. Since $|\mathcal I|$ is even, choose a partition
$\mathcal P$ of $\mathcal I$ into two-element sets.

For each $P=\{(u,i),(v,j)\}\in\mathcal P$, add a copy of
$K_{D+1}$ to $G_0$. Choose two distinct vertices $x_P,y_P$ in
this copy, delete the edge $x_Py_P$, and add the edges $ux_P$
and $vy_P$. Denote the resulting graph by $\widehat G$.

Each index $(v,i)$ contributes exactly one edge incident with
$v$. Thus every $v\in V(G_0)$ has degree
\[
 d_{\widehat G}(v)=d_{G_0}(v)+|\{i:(v,i)\in\mathcal I\}|=D.
\]
For each $P$, deleting $x_Py_P$ reduces the degrees of $x_P$ and
$y_P$ to $D-1$, and the two edges to $G_0$ restore both degrees to
$D$. All other vertices of the added copies have degree $D$.
Hence $\widehat G$ is $D$-regular.

Each added copy is simple, and each of its vertices is incident
with at most one edge to $G_0$. Thus $\widehat G$ has no loops
and every edge added to $G$ has multiplicity one. All added edges
have an endpoint outside $V(G)$, so $\widehat G[V(G)]=G$.
\end{proof}

Apply \textbf{Lemma~\ref{lem:regularization}} to the reduced
multigraph. The induced embedding preserves the marked strong clique.
We keep the original marks, leave every added edge unmarked, and
write $G$ for the completed multigraph.
Let $G_{\mathrm s}$ be the underlying simple graph of $G$. Color
an edge $uv$ of $G_{\mathrm s}$ red if a marked edge of $G$ has
endpoints $u,v$, and blue otherwise. Let $s$ be the number of red
edges. For each red edge $uv$, choose one marked edge with endpoints
$u,v$ as its \emph{representative}. Let $b_{uv}$ be the number of
remaining marked edges with endpoints $u,v$. In total there are
$b=m-s$ remaining marked edges. We use the quantities
\begin{equation}\label{eq:weighted-parameters}
 M=\frac{m}{D^2},\qquad S=\frac{s}{D^2},\qquad
 e=\frac{b}{m},\qquad S=M(1-e).
\end{equation}

For each integer $n\geq0$, let $P_n=\{0,\ldots,n-1\}$, with
$P_0=\varnothing$. We call the elements of $P_n$ \emph{positions}.
A \emph{pattern} $F$ on $P_n$ assigns a color
$c_F(r,t)\in\{0,1,2\}$ to each unordered pair $\{r,t\}$ of distinct
positions. Its \emph{underlying graph} has vertex set $P_n$ and an edge
between $r$ and $t$ whenever $c_F(r,t)\in\{1,2\}$. We call an edge of
color $2$ in the underlying graph of $F$ red. Given a pair of vertices
$u,v\in V(G_{\mathrm s})$, define $c_G(u,v)$ to be $2$ if $uv$ is red,
$1$ if $uv$ is blue, and $0$ if $u=v$ or $u,v$ are nonadjacent. A map
$\phi\colon P_n\to V(G_{\mathrm s})$ \emph{gives the colors of $F$} if
$c_G(\phi(r),\phi(t))=c_F(r,t)$ for every position pair $\{r,t\}$.
Note that such a map may send distinct positions to the same vertex.

Suppose $\phi$ gives the colors of $F$ and that $\{a,b\}$ and
$\{c,d\}$ are disjoint red edges of $F$. There is a marked edge with
endpoints $\phi(a),\phi(b)$ and another with endpoints
$\phi(c),\phi(d)$. If their endpoint sets are disjoint, the strong
clique property gives an edge of $G$ joining the sets. Its endpoints
are adjacent in $G_{\mathrm s}$. If the sets meet, exchange the
positions within each pair if needed so that $\phi(a)=\phi(c)$.
Then $\phi(a)\phi(d)$ is the same red edge of $G_{\mathrm s}$ as
$\phi(c)\phi(d)$. Thus
\begin{equation}\label{eq:weighted-joining-edge}
 \max_{\substack{r\in\{a,b\}\\t\in\{c,d\}}}c_F(r,t)\geq1.
\end{equation}
In particular, the red edges of a pattern given by a map form a
strong clique in its underlying graph.

A \emph{weighted pattern} on $P_n$ is a pair $(F,J)$, where $F$ is a
pattern on $P_n$ and $J$ is a set of red unordered pairs of positions
such that the graph on $P_n$ with edge set $J$ is a forest. We call the
elements of $J$ the \emph{specified pairs}. For a weighted pattern
$(F,J)$ on $P_n$, define its \emph{weighted density} by
\begin{equation}\label{eq:weighted-density}
 Q(F,J)=\frac{1}{2mD^{n-2}}\sum_{\phi}
             \prod_{\{r,t\}\in J} b_{\phi(r)\phi(t)},
\end{equation}
where the sum runs over the maps giving the colors of $F$. An empty
product is $1$. If $\ell$ specified pairs have the same unordered
image $\{u,v\}$, their factors in the product give $b_{uv}^{\ell}$.
Indeed, for a fixed map $\phi$, the product counts the choices of
one remaining marked edge with endpoints $\phi(r),\phi(t)$ for
each specified pair $\{r,t\}$. The choices are independent, and the
same edge may be chosen for different pairs. Let $h(F,J)$ count the
tuples $\bigl(\phi,(\eta_p)_{p\in J}\bigr)$ consisting of such a map
and its chosen edges $\eta_p$. Thus
\begin{equation}\label{eq:weighted-raw-count}
 h(F,J)=2mD^{n-2}Q(F,J).
\end{equation}

Let $R=\{\{0,1\},\{2,3\}\}$ and let $\mathcal B$ consist of all
patterns $F$ on $P_4$ for which the pairs in $R$ are red and the
red edges form a strong clique in the underlying graph. The
underlying graph of each $F\in\mathcal B$ is connected because
the two red pairs in $R$ cover $P_4$ and an edge joins them. Define
\[
 \mathcal M(Q)=\frac{1}{2}\sum_{F\in\mathcal B}
                           \sum_{J\subseteq R} Q(F,J).
\]

For a loopless multigraph of maximum degree $D$,
\textbf{Conjecture~\ref{conj:fgst}} asks us to show that a strong
clique of size $m$ satisfies $m\leq f(D)$, or equivalently,
\[
 M=\frac{m}{D^2}\leq\frac{f(D)}{D^2}.
\]
The following theorem expresses $M$ in terms of the weighted
densities as the linear combination $\mathcal M(Q)$. This reduces
the problem of bounding $m$ to that of bounding $\mathcal M(Q)$.

\begin{theorem}\label{thm:weighted-objective}
Fix a positive integer $D$. Given a nonempty marked strong clique
in a finite loopless multigraph of maximum degree at most $D$,
compute the weighted densities after the reduction and regular
completion above. Then
\begin{equation}\label{eq:weighted-objective}
 M=\mathcal M(Q).
\end{equation}
The following statements are equivalent.
\begin{enumerate}[label=\textup{(\roman*)}]
 \item Every strong clique in every finite loopless multigraph of
       maximum degree at most $D$ has at most $f(D)$ edges.
 \item Every collection of weighted densities obtained from a
       nonempty marked strong clique with degree bound $D$
       satisfies $\mathcal M(Q)\leq\frac{f(D)}{D^2}$.
 \item No collection of weighted densities obtained from a
       nonempty marked strong clique with degree bound $D$
       satisfies $\mathcal M(Q)=\frac{f(D)+1}{D^2}$.
\end{enumerate}
\end{theorem}

\begin{proof}
We count pairs of marked edges with an endpoint order
for each edge. For each $q=\{r,t\}\in R$ with $r<t$, choose a
marked edge $\varepsilon_q$ and assign its endpoints to $r$ and
$t$ in either order. Each pair in $R$ permits $2m$ choices.
The choices are independent and may use the same edge, so there
are $(2m)^2$ choices in total.

The assignments define a map $\phi\colon P_4\to V(G_{\mathrm s})$.
Define $F$ by $c_F(r,t)=c_G(\phi(r),\phi(t))$ for every pair
$\{r,t\}$ in $P_4$. Each pair in $R$ is red because its chosen
edge is marked. By \eqref{eq:weighted-joining-edge}, the red edges
of $F$ form a strong clique in its underlying graph, so
$F\in\mathcal B$. Set
\[
 J=\{q\in R:\varepsilon_q\text{ is a remaining marked edge}\},
 \qquad \eta_q=\varepsilon_q\quad(q\in J).
\]
Every subset of $R$ is the edge set of a forest, since the pairs in
$R$ are disjoint.
Thus $\bigl(\phi,(\eta_q)_{q\in J}\bigr)$ is counted by $h(F,J)$.

Conversely, take $F\in\mathcal B$, $J\subseteq R$ and a tuple
$\bigl(\phi,(\eta_q)_{q\in J}\bigr)$ counted by $h(F,J)$. For
each $q=\{r,t\}\in R$ with $r<t$, define $\varepsilon_q$ to be
$\eta_q$ if $q\in J$ and the representative with endpoints
$\phi(r),\phi(t)$ if $q\in R\setminus J$. The representative
exists because $q$ is red in $F$. In both cases take
$(\phi(r),\phi(t))$ as the endpoint order. The two constructions
are inverse, since $(F,J,\phi,(\eta_q)_{q\in J})$ determines
each $\varepsilon_q$ and its endpoint order, and the resulting
choices give the same $F$, $J$, $\phi$ and $\eta_q$. Thus
\eqref{eq:weighted-raw-count} gives
\[
 (2m)^2=\sum_{F\in\mathcal B}\sum_{J\subseteq R}h(F,J)
       =2mD^2\sum_{F\in\mathcal B}\sum_{J\subseteq R}Q(F,J).
\]
Dividing by $4mD^2$ gives \eqref{eq:weighted-objective}.

If \textup{(i)} holds, every nonempty marked strong clique used
in \textup{(ii)} satisfies $m\leq f(D)$. Hence
\eqref{eq:weighted-objective} gives \textup{(ii)}. Statement
\textup{(ii)} implies \textup{(iii)} because
$\frac{f(D)}{D^2}<\frac{f(D)+1}{D^2}$.
To prove that \textup{(iii)} implies \textup{(i)}, suppose that a
strong clique has more than $f(D)$ edges. Mark exactly $f(D)+1$ of
its edges and compute the weighted densities with degree bound $D$.
The marked edges form a strong clique, and
\eqref{eq:weighted-objective} gives
$\mathcal M(Q)=\frac{f(D)+1}{D^2}$, contrary to \textup{(iii)}.
\end{proof}

To bound $\mathcal M(Q)$, we first derive two basic constraints
on the weighted densities.

Let $F_\bullet$ be the pattern on $P_2$ whose only position pair
$p_\bullet=\{0,1\}$ is red. To specify a map $\phi$ giving the colors
of $F_\bullet$, we give the ordered pair $(\phi(0),\phi(1))$. Its
first entry is the image of position $0$, and its second entry is the
image of position $1$. Thus each red edge $uv$ of $G_{\mathrm s}$
gives two maps, one with $(\phi(0),\phi(1))=(u,v)$ and one with
$(\phi(0),\phi(1))=(v,u)$.

When $J=\varnothing$, each map contributes one choice, so
$h(F_\bullet,\varnothing)=2s$. When $J=\{p_\bullet\}$, a choice
consists of a map $\phi$ together with a remaining marked edge whose
endpoints are $\phi(0)$ and $\phi(1)$. Choose one of the $b$ remaining
marked edges of $G$, and assign its endpoints to positions $0$ and $1$
in either of the two ways above. Thus $h(F_\bullet,\{p_\bullet\})=2b$.
Since $F_\bullet$ has two positions, both weighted densities are
obtained by dividing the corresponding counts by $2m$. Using
$s+b=m$, we obtain
\begin{equation}\label{eq:weighted-normalization}
 Q(F_\bullet,\varnothing)=\frac{s}{m},\qquad
 Q(F_\bullet,\{p_\bullet\})=\frac{b}{m},\qquad
 Q(F_\bullet,\varnothing)+Q(F_\bullet,\{p_\bullet\})=1.
\end{equation}

\begin{lemma}\label{lem:weighted-density-bounds}
The weighted density $Q(F,J)$ is nonnegative. If the underlying
graph of $F$ is connected and has a red edge, then
\begin{equation}\label{eq:weighted-density-bounds}
 Q(F,J)\leq
 \begin{cases}
  1-e,& J=\varnothing,\\
  e,& J\ne\varnothing .
 \end{cases}
\end{equation}
\end{lemma}

\begin{proof}
Nonnegativity follows from \eqref{eq:weighted-density}, since
every $b_{uv}$ is nonnegative. To prove the upper bounds, write
$P_n$ for the position set of $F$. If $J=\varnothing$, choose a red
pair $p$ of $F$. If $J\ne\varnothing$, choose $p\in J$. Let $U$ be
the forest on $P_n$ with edge set $J\cup\{p\}$.

To bound $h(F,J)$, we choose a map $\phi\colon P_n\to V(G_{\mathrm s})$
and a remaining marked edge for each pair in $J$. We choose the
images under $\phi$ in the following order. Begin with the two
positions in $p$ and set $A_2=p$. For $i=2,\ldots,n-1$, if $U$ has
an edge between $A_i$ and $P_n\setminus A_i$, choose such an edge
as $q_i$. Otherwise, choose $q_i$ to be an edge of the underlying
graph of $F$ between these two sets. The latter choice is possible
by connectedness. Write $q_i=\{r_i,t_i\}$ with $r_i\in A_i$ and
$t_i\in P_n\setminus A_i$, and set $A_{i+1}=A_i\cup\{t_i\}$.
After choosing the images of the positions in $A_2$, we choose
$\phi(t_2),\ldots,\phi(t_{n-1})$ in that order.

We show that $p,q_2,\ldots,q_{n-1}$ include all edges of $U$.
If $q_i\notin E(U)$, no edge of $U$ joins $A_i$ to its complement.
Hence $A_i$ is a union of vertex sets of components of $U$, and
$t_i$ is the first position chosen in its component. Conversely,
if $t_i$ is the first position chosen in its component, then
$r_i\in A_i$ belongs to another component, so $q_i\notin E(U)$.
If $U$ has $c$ components, exactly $c-1$ of the edges $q_i$ lie
outside $E(U)$, one for each component other than the one containing $p$.
The edges $p,q_2,\ldots,q_{n-1}$ are distinct, since each $q_i$ contains
the new position $t_i$. Exactly $n-c$ of them lie in $E(U)$, which
is the number of edges in the forest $U$. In particular, every pair in
$J$ is either $p$ or one of the $q_i$.

To choose the images of the positions in $p$, take a red edge of
$G_{\mathrm s}$ if $J=\varnothing$, or a remaining marked edge of $G$
if $J\ne\varnothing$. Assign the endpoints to the two positions in
either order. This gives $2s$ or $2b$ choices, respectively.

At step $i$, the image $\phi(r_i)$ has already been chosen. If
$q_i\in J$, choose a remaining marked edge of $G$ incident with
$\phi(r_i)$ and set $\phi(t_i)$ to its other endpoint. If
$q_i\notin J$, choose $\phi(t_i)$ among the neighbors of $\phi(r_i)$
in $G_{\mathrm s}$. Each step has at most $D$ choices because $G$
is $D$-regular and $G_{\mathrm s}$ has maximum degree at most $D$.

For a tuple $\bigl(\phi,(\eta_q)_{q\in J}\bigr)$ counted by $h(F,J)$,
use $\phi$ to choose the images and $\eta_q$ to choose the edge for
each pair $q\in J$. These choices are allowed because $\phi$ gives
the colors of $F$ and each $\eta_q$ is a remaining marked edge with
the prescribed endpoints. The choices determine the tuple because
they specify every vertex image and every edge $\eta_q$.
Since there are $n-2$ steps after the initial choice, we obtain
\[
 h(F,\varnothing)\leq2sD^{n-2},\qquad
 h(F,J)\leq2bD^{n-2}\quad\text{if }J\ne\varnothing .
\]
Divide by $2mD^{n-2}$ and use \eqref{eq:weighted-raw-count},
$\frac{s}{m}=1-e$ and $\frac{b}{m}=e$ to obtain
\eqref{eq:weighted-density-bounds}.
\end{proof}

The two basic constraints above are insufficient to obtain the bound
in \textbf{Conjecture~\ref{conj:fgst}} from
\textbf{Theorem~\ref{thm:weighted-objective}}. To obtain sharper
estimates of $\mathcal M(Q)$, we introduce two tools in
\textbf{Subsections~\ref{subsec:weighted-flag-inequalities}
and~\ref{subsec:weighted-density-constraints}}, namely
\emph{weighted flag inequalities} and \emph{weighted density constraints}.

By \textbf{Theorem~\ref{thm:weighted-objective}}, we can prove
\textbf{Conjecture~\ref{conj:fgst}} by establishing
$\mathcal M(Q)\leq\frac{f(D)}{D^2}$. We may estimate $\mathcal M(Q)$
directly or suppose that it lies in $\bigl(\frac{f(D)}{D^2},\infty\bigr)$,
narrow this range using preliminary estimates, and partition the
remaining range into intervals. On each interval, we use the assumed
bounds on $M=\mathcal M(Q)$ to analyze the structure of $G$. We turn
the resulting structural restrictions into weighted density constraints
for that interval and combine them with the weighted flag inequalities.
An upper bound on $\mathcal M(Q)$ below the left endpoint rules out
the interval.

For odd $D$, we use part \textup{(iii)} of
\textbf{Theorem~\ref{thm:weighted-objective}} and assume that there is
a strong clique with exactly $m=f(D)+1$ edges. We impose
\[
 \mathcal M(Q)=\frac{f(D)+1}{D^2}
\]
as a weighted density constraint. We rule out the assumed clique by
proving an upper bound on $\mathcal M(Q)$ below this value or by
showing that the weighted flag inequalities and the weighted density
constraints for the given $D$ have no common solution.

\subsection{Weighted Flag Inequalities}\label{subsec:weighted-flag-inequalities}

We introduce weighted types and weighted flags and use them to
derive weighted flag inequalities for the weighted densities.
Fix a weighted pattern $(T,J_0)$ on $P_k$. A weighted pattern
$(F,J)$ on $P_n$ is a \emph{weighted flag of $(T,J_0)$} if
$n\geq k$ and
\[
 F|_{P_k}=T,\qquad J\cap\binom{P_k}{2}=J_0,
\]
where $F|_{P_k}$ restricts the colors of $F$ to pairs in $P_k$.
In this case, we call $(T,J_0)$ the \emph{weighted type} of the flag
and say that $(F,J)$ is a weighted flag of type $(T,J_0)$.

Fix a weighted pattern $(T,J_0)$ on $P_k$ and a list of weighted
flags $(F_1,J_1),\ldots,(F_\ell,J_\ell)$ of weighted type $(T,J_0)$,
where $(F_i,J_i)$ is on $P_{n_i}$.
For each map $\phi\colon P_k\to V(G_{\mathrm s})$ giving the
colors of $T$, define
\[
 w(\phi)=\prod_{\{r,t\}\in J_0} b_{\phi(r)\phi(t)}.
\]
This product counts the choices of a remaining marked edge
with endpoints $\phi(r),\phi(t)$ for each pair $\{r,t\}\in J_0$.
An \emph{extension of $\phi$ to $F_i$} is a map
$g_i\colon P_{n_i}\to V(G_{\mathrm s})$ such that
$g_i|_{P_k}=\phi$ and $g_i$ gives the colors of $F_i$. Define
\[
 a_i(\phi)=\sum_{g_i}
     \prod_{\{r,t\}\in J_i\setminus J_0} b_{g_i(r)g_i(t)},
\]
where the sum ranges over these extensions. For each extension,
the product counts the choices of a remaining marked edge for
each specified pair in $J_i\setminus J_0$, with endpoints given
by $g_i$. Let $v(\phi)$ be the column vector with entries
$v_i(\phi)=\frac{a_i(\phi)}{D^{n_i-k}}$. The \emph{weighted matrix}
of the list of flags is
\begin{subequations}
\begin{equation}\label{eq:weighted-matrix}
 W=\frac{1}{2mD^{k-2}}
       \sum_{\phi} w(\phi)v(\phi)v(\phi)^{\mathsf T},
\end{equation}
where the sum runs over the maps giving the colors of $T$.

Note that $W$ is positive semidefinite. Since $w(\phi)\geq0$ for
every $\phi$, for every $x\in\mathbb R^\ell$ we have
\begin{equation}\label{eq:weighted-matrix-psd}
 x^{\mathsf T}Wx
 =\frac{1}{2mD^{k-2}}\sum_\phi
      w(\phi)\bigl(x^{\mathsf T}v(\phi)\bigr)^2\geq0.
\end{equation}
\end{subequations}

To describe the entries of $W$, fix $i,j$ and put
\[
 n=n_i+n_j-k,\qquad A=P_{n_i},\qquad
 B=P_k\cup\{n_i,\ldots,n-1\}.
\]
Keep the positions of $(F_i,J_i)$ in $A$ and relabel the positions
of $(F_j,J_j)$ by the bijection $\rho\colon P_{n_j}\to B$ given by
\begin{equation}\label{eq:flag-relabeling}
 \rho(t)=
 \begin{cases}
  t,& t\in P_k,\\
  n_i+t-k,& t\in P_{n_j}\setminus P_k .
 \end{cases}
\end{equation}
The relabeled pattern assigns the color $c_{F_j}(r,t)$ to
$\{\rho(r),\rho(t)\}$. Its specified pairs are
$J_j^B=\{\{\rho(r),\rho(t)\}:\{r,t\}\in J_j\}$.
Write $J_{ij}=J_i\cup J_j^B$. Since $A\cap B=P_k$ and the specified
pairs of each flag within $P_k$ are exactly those in $J_0$, we
have $J_i\cap J_j^B=J_0$. Hence
\[
 J_{ij}=J_0\ \dot\cup\ (J_i\setminus J_0)\
              \dot\cup\ (J_j^B\setminus J_0).
\]

Let $\mathcal U_{ij}$ consist of the patterns $H$ on $P_n$ that
agree with $F_i$ on $A$ and with the relabeled $F_j$ on $B$.
For each pair with one position in $A\setminus P_k$ and the other
in $B\setminus P_k$, allow any color in $\{0,1,2\}$. Retain the
patterns whose red edges form a strong clique in the underlying
graph, with the condition holding vacuously when there are no red
edges. Each specified pair in $J_{ij}$ is red in every remaining
pattern, since it belongs to one of the two flags. If $J_{ij}$ is
a forest, then every $(H,J_{ij})$ with $H\in\mathcal U_{ij}$ is a
weighted pattern.

\begin{theorem}\label{thm:weighted-flag-inequality}
If $J_{ij}$ is a forest, then
\begin{equation}\label{eq:weighted-matrix-entry}
 W_{ij}=\sum_{H\in\mathcal U_{ij}}Q(H,J_{ij}).
\end{equation}
\end{theorem}

\begin{proof}
Fix a map $\phi$ giving the colors of $T$. In the sum for $a_i(\phi)$,
write the extension as $g_i$ on $A$. In the sum for $a_j(\phi)$,
write $g_j$ for the map on $B$ obtained by relabeling an extension.
The corresponding extension on $P_{n_j}$ sends $t$ to
$g_j(\rho(t))$. Both maps equal $\phi$ on $A\cap B=P_k$, so they determine a
unique map $g\colon P_n\to V(G_{\mathrm s})$ with $g|_A=g_i$ and
$g|_B=g_j$. Using the disjoint union of specified pairs defining
$J_{ij}$, we expand $w(\phi)a_i(\phi)a_j(\phi)$ as
\begin{equation}\label{eq:weighted-flag-fixed-type}
\begin{aligned}
 w(\phi)a_i(\phi)a_j(\phi)
 &=\sum_{g_i,g_j}
   \left(\prod_{\{r,t\}\in J_0} b_{\phi(r)\phi(t)}\right)
   \left(\prod_{\{r,t\}\in J_i\setminus J_0} b_{g_i(r)g_i(t)}\right)\\
 &\qquad{}\times
   \left(\prod_{\{r,t\}\in J_j^B\setminus J_0} b_{g_j(r)g_j(t)}\right)\\
 &=\sum_{g_i,g_j}\prod_{\{r,t\}\in J_{ij}}b_{g(r)g(t)}.
\end{aligned}
\end{equation}

For each choice of $\phi,g_i,g_j$ in \eqref{eq:weighted-flag-fixed-type},
the combined map $g$ determines a pattern $H$ on $P_n$. The
restrictions of $g$ to $A$ and $B$ give the colors prescribed by $F_i$
and the relabeled $F_j$, respectively. By \eqref{eq:weighted-joining-edge},
the red edges of $H$ form a strong clique in its underlying graph.
Thus $H\in\mathcal U_{ij}$.

Conversely, let $H\in\mathcal U_{ij}$ and let $g$ give the colors of
$H$. Set $\phi=g|_{P_k}$, $g_i=g|_A$, and $g_j=g|_B$. Since $H$
agrees with $F_i$ on $A$ and with the relabeled $F_j$ on $B$, $g_i$
and $g_j$ give the colors of these two patterns, respectively.
Moreover, $g_i|_{P_k}=g_j|_{P_k}=\phi$, and $\phi$ gives the colors
of $T$ because $F_i|_{P_k}=T$. These restrictions uniquely recover
the three maps that produce $g$. We can therefore first sum
\eqref{eq:weighted-flag-fixed-type} over $\phi$, then group the terms
by the pattern $H$ determined by $g$. This gives
\begin{equation}\label{eq:weighted-flag-grouping}
\begin{aligned}
 \sum_\phi w(\phi)a_i(\phi)a_j(\phi)
 &=\sum_\phi\sum_{g_i,g_j}
      \prod_{\{r,t\}\in J_{ij}}b_{g(r)g(t)}\\
 &=\sum_{H\in\mathcal U_{ij}}
    \sum_{g\text{ giving the colors of }H}
      \prod_{\{r,t\}\in J_{ij}}b_{g(r)g(t)}.
\end{aligned}
\end{equation}

When $J_{ij}$ is a forest, $(H,J_{ij})$ is a weighted pattern for
every $H\in\mathcal U_{ij}$. Substituting the definitions of
$v_i(\phi)$ and $v_j(\phi)$ into \eqref{eq:weighted-matrix}, and
using \eqref{eq:weighted-flag-grouping} and \eqref{eq:weighted-density},
we obtain
\[
\begin{aligned}
 W_{ij}
 &=\frac{1}{2mD^{k-2}}\sum_\phi
   w(\phi)\frac{a_i(\phi)}{D^{n_i-k}}\frac{a_j(\phi)}{D^{n_j-k}}\\
 &=\frac{1}{2mD^{n-2}}\sum_\phi w(\phi)a_i(\phi)a_j(\phi)\\
 &=\sum_{H\in\mathcal U_{ij}}\frac{1}{2mD^{n-2}}
    \sum_{g\text{ giving the colors of }H}
      \prod_{\{r,t\}\in J_{ij}}b_{g(r)g(t)}\\
 &=\sum_{H\in\mathcal U_{ij}}Q(H,J_{ij}),
\end{aligned}
\]
where the second equality uses $(k-2)+(n_i-k)+(n_j-k)=n-2$.
This proves \eqref{eq:weighted-matrix-entry}.
\end{proof}

The forest condition in \eqref{eq:weighted-matrix-entry} can fail.
For example, take $T$ to be the red pair on $P_2$ and let
$J_0=\varnothing$. Use two red triangles as the flags and specify
the two pairs joining the third position of each triangle to $0$
and $1$. Under \eqref{eq:flag-relabeling}, the third positions become
$2$ and $3$. The specified pairs in $J_{ij}$ then form the cycle
$0,2,1,3,0$.

When every $J_{ij}$ is a forest, substituting
\eqref{eq:weighted-matrix-entry} into the quadratic form in
\eqref{eq:weighted-matrix-psd} gives, for each fixed
$x\in\mathbb R^\ell$, the linear inequality
\[
 \sum_{i,j=1}^{\ell}x_ix_j
      \sum_{H\in\mathcal U_{ij}}Q(H,J_{ij})\geq0 .
\]
We call each such inequality a \emph{weighted flag inequality}.

\subsection{Weighted Density Constraints}\label{subsec:weighted-density-constraints}

A \emph{weighted density constraint} is a linear equation or
inequality in the weighted densities $Q(F,J)$. Some weighted density
constraints hold for every value of $M$. Others follow from structural
properties of $G$ that can be deduced under an assumption on $M$.
As examples, we first derive two constraints of the first kind.

For a pattern $F$ on $P_k$, a position $u\in P_k$, and $c\in\{1,2\}$,
let $\mathcal E_c(F,u)$ consist of the patterns $H$ on $P_{k+1}$
such that $H|_{P_k}=F$, $c_H(u,k)=c$, and the red edges of $H$ form
a strong clique in its underlying graph.

\begin{theorem}\label{thm:weighted-extension}
Let $(F,J)$ be a weighted pattern on $P_k$ whose underlying
graph is connected and in which the red edges form a nonempty
strong clique. Fix $u\in P_k$ and put $p=\{u,k\}$. Then
\begin{equation}\label{eq:weighted-extension}
 \sum_{H\in\mathcal E_1(F,u)}Q(H,J)
 +\sum_{H\in\mathcal E_2(F,u)}
    \bigl(Q(H,J)+Q(H,J\cup\{p\})\bigr)=Q(F,J).
\end{equation}
\end{theorem}

\begin{proof}
An \emph{extension of $(F,J)$ at $u$} consists of a tuple
$\bigl(\phi,(\eta_q)_{q\in J},\varepsilon\bigr)$, where
$\bigl(\phi,(\eta_q)_{q\in J}\bigr)$ is counted by $h(F,J)$ and
$\varepsilon$ is an edge of $G$ incident with $\phi(u)$.
There are $D\,h(F,J)$ extensions, since $G$ is $D$-regular and
parallel edges count as distinct choices of $\varepsilon$.

To count the extensions by their patterns on $P_{k+1}$, define
$g\colon P_{k+1}\to V(G_{\mathrm s})$ by $g|_{P_k}=\phi$, and
let $g(k)$ be the other endpoint of $\varepsilon$. Let $H$ be
the pattern given by $g$, and put $c=c_H(u,k)$. We have
$H|_{P_k}=F$, and \eqref{eq:weighted-joining-edge} implies that the
red edges of $H$ form a strong clique in its underlying graph.
Since $g(u)$ and $g(k)$ are joined by $\varepsilon$, we have
$c\in\{1,2\}$ and hence $H\in\mathcal E_c(F,u)$.

An unmarked $\varepsilon$ gives $c=1$, since the reduction removes
every unmarked edge parallel to a marked edge and the regular
completion adds no parallel edges. A marked $\varepsilon$ gives
$c=2$. Define the set of specified pairs $J'$ by the following cases.
\[
 (c,J')=
 \begin{cases}
 (1,J),&\text{if $\varepsilon$ is unmarked},\\
 (2,J),&\text{if $\varepsilon$ is a representative},\\
 (2,J\cup\{p\}),&\text{if $\varepsilon$ is a remaining marked edge}.
 \end{cases}
\]
Take $\xi_q=\eta_q$ for $q\in J$ and, if $p\in J'$, take
$\xi_p=\varepsilon$. When $p$ is added, it is red because $c=2$,
and it is the only edge of $J'$ incident with the new position $k$.
Thus $J'$ is the edge set of a forest. The resulting tuple
$\bigl(g,(\xi_q)_{q\in J'}\bigr)$ is counted by $h(H,J')$.
\textbf{Figure~\ref{fig:weighted-extension-choices}} shows the
three possibilities.

\begin{figure}[htbp]
\centering
\begin{tikzpicture}[every node/.style={font=\small},
                    marked edge/.append style={draw=red!70!black},
                    unmarked edge/.append style={draw=blue!65!black,dashed}]
  \node at (1.03,1.31) {Unmarked edge};
  \node[graph vertex,label=below:$\phi(u)$] (a) at (0.24,0.46) {};
  \node[graph vertex,label=below:$g(k)$] (b) at (1.82,0.46) {};
  \draw[unmarked edge] (a)--node[above] {$\varepsilon$} (b);
  \node at (1.03,-0.37) {$(c_H(u,k),J')=(1,J)$};
  \begin{scope}[xshift=4.55cm]
    \node at (1.03,1.31) {Representative};
    \node[graph vertex,label=below:$\phi(u)$] (a) at (0.24,0.46) {};
    \node[graph vertex,label=below:$g(k)$] (b) at (1.82,0.46) {};
    \draw[marked edge] (a)--node[above] {$\varepsilon$} (b);
    \node at (1.03,-0.37) {$(c_H(u,k),J')=(2,J)$};
  \end{scope}
  \begin{scope}[xshift=9.25cm]
    \node at (1.03,1.31) {Remaining marked edge};
    \node[graph vertex,label=below:$\phi(u)$] (a) at (0.24,0.46) {};
    \node[graph vertex,label=below:$g(k)$] (b) at (1.82,0.46) {};
    \draw[marked edge] (a)--(b);
    \draw[marked edge] (a) to[bend left=22]
     node[above] {$\varepsilon$} (b);
    \node at (1.03,-0.37) {$(c_H(u,k),J')=(2,J\cup\{p\})$};
  \end{scope}
\end{tikzpicture}
\caption{An extension along the chosen edge $\varepsilon$. Here
$J'$ is the set of specified pairs after the extension. In the last panel,
the straight edge is the representative.}
\label{fig:weighted-extension-choices}
\end{figure}
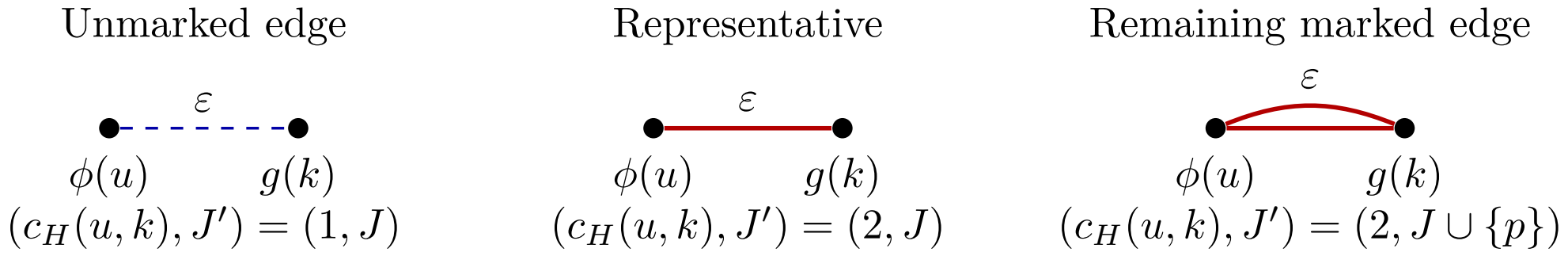

Conversely, fix $H\in\mathcal E_c(F,u)$ and let $J'$ be one of the
sets in the cases above for this value of $c$. Given a tuple
$\bigl(g,(\xi_q)_{q\in J'}\bigr)$ counted by $h(H,J')$, set
$\phi=g|_{P_k}$ and $\eta_q=\xi_q$ for $q\in J$. The edge
$\varepsilon$ must have endpoints $g(u),g(k)$. If $(c,J')=(1,J)$,
take the unique unmarked edge with these endpoints. Its uniqueness
follows from the reduction and from the fact that every added edge
in the regular completion has multiplicity one. If $(c,J')=(2,J)$,
take the representative with these endpoints. If
$(c,J')=(2,J\cup\{p\})$, take $\varepsilon=\xi_p$. In each case
the resulting extension gives the original $g$ and $(\xi_q)_{q\in J'}$.
The extension is unique, since the restrictions give $\phi$ and
$(\eta_q)_{q\in J}$ and the case determines $\varepsilon$.
Thus exactly $h(H,J')$ extensions give $H$ and $J'$.

Summing over $H$ and the corresponding choices of $J'$ counts all
$D\,h(F,J)$ extensions. Hence
\[
 \sum_{H\in\mathcal E_1(F,u)}h(H,J)
 +\sum_{H\in\mathcal E_2(F,u)}
     \bigl(h(H,J)+h(H,J\cup\{p\})\bigr)=D\,h(F,J).
\]
Dividing by $2mD^{k-1}$ and applying \eqref{eq:weighted-raw-count}
gives \eqref{eq:weighted-extension}.
\end{proof}

\begin{theorem}\label{thm:remaining-edge-orderings}
With $R$ and $\mathcal B$ as in
\textbf{Subsection~\ref{subsec:weighted-sketch}}, let
$\mathcal B_+=\{F\in\mathcal B:c_F(1,2)\in\{1,2\}\}$.
Put $J_\star=\{\{0,1\},\{0,2\}\}$. For $c\in\{0,1,2\}$ let
$H_c$ be the pattern on $P_3$ whose pairs in $J_\star$ have
color $2$ and whose pair $\{1,2\}$ has color $c$. Then
\begin{equation}\label{eq:remaining-edge-orderings}
 4\sum_{F\in\mathcal B_+}Q(F,R)
 \geq\sum_{F\in\mathcal B}Q(F,R)
     +\frac{2}{D}\sum_{c=0}^2 Q(H_c,J_\star).
\end{equation}
\end{theorem}

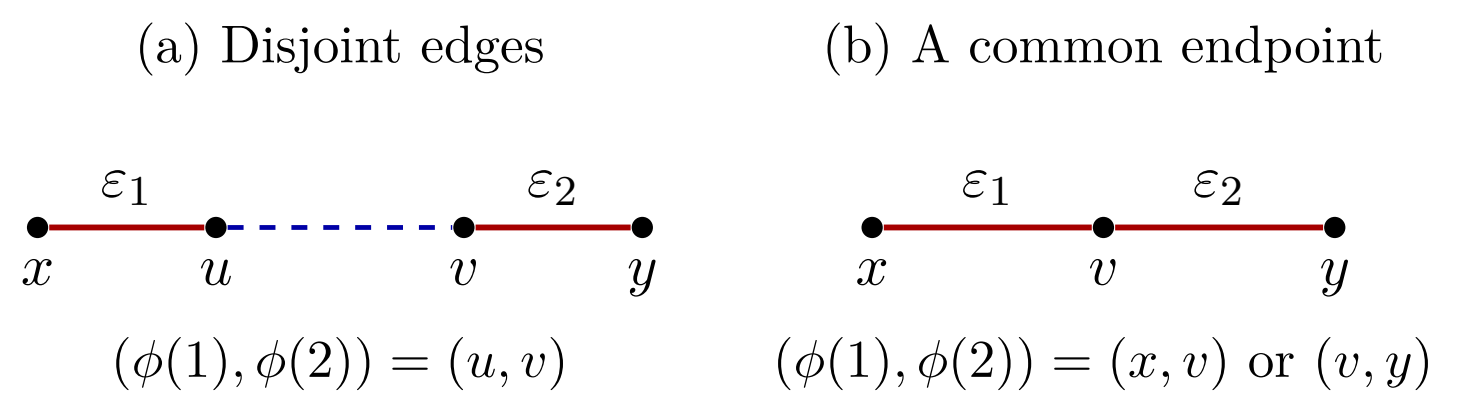
\begin{figure}[htbp]
\centering
\begin{tikzpicture}[
  vertex/.style={circle,fill=black,inner sep=0pt,minimum size=4pt},
  selected edge/.style={draw=red!65!black,line width=1.1pt},
  joining edge/.style={draw=blue!65!black,dashed,line width=0.9pt},
  panel title/.style={font=\small},
  order label/.style={font=\small}
]
\node[panel title] at (2.05,1.18) {(a) Disjoint edges};
\node[vertex,label=below:$x$] (disjoint-x) at (0.05,0) {};
\node[vertex,label=below:$u$] (disjoint-u) at (1.23,0) {};
\node[vertex,label=below:$v$] (disjoint-v) at (2.87,0) {};
\node[vertex,label=below:$y$] (disjoint-y) at (4.05,0) {};
\draw[selected edge] (disjoint-x) -- node[above] {$\varepsilon_1$} (disjoint-u);
\draw[joining edge] (disjoint-u) -- (disjoint-v);
\draw[selected edge] (disjoint-v) -- node[above] {$\varepsilon_2$} (disjoint-y);
\node[order label] at (2.05,-0.9) {$(\phi(1),\phi(2))=(u,v)$};

\begin{scope}[xshift=5.05cm]
  \node[panel title] at (2.05,1.18) {(b) A common endpoint};
  \node[vertex,label=below:$x$] (shared-x) at (0.52,0) {};
  \node[vertex,label=below:$v$] (shared-v) at (2.05,0) {};
  \node[vertex,label=below:$y$] (shared-y) at (3.58,0) {};
  \draw[selected edge] (shared-x) -- node[above] {$\varepsilon_1$} (shared-v);
  \draw[selected edge] (shared-v) -- node[above] {$\varepsilon_2$} (shared-y);
  \node[order label] at (2.05,-0.9)
    {$(\phi(1),\phi(2))=(x,v)\text{ or }(v,y)$};
\end{scope}
\end{tikzpicture}
\caption{Endpoint orders for two remaining marked edges. In (a), the
dashed edge joins the endpoints mapped from positions $1$ and $2$.
In (b), the common endpoint may be the image of either of those
positions.}
\label{fig:remaining-edge-orderings}
\end{figure}

\begin{proof}
Choose an ordered pair $(\varepsilon_1,\varepsilon_2)$ of
remaining marked edges, allowing $\varepsilon_1=\varepsilon_2$.
There are $b^2$ such pairs. For each pair, choose an endpoint
order for $\varepsilon_1$ at positions $0,1$ and one for
$\varepsilon_2$ at positions $2,3$. Let $\phi\colon P_4\to
V(G_{\mathrm s})$ be the map given by these endpoint assignments.
By the counting correspondence
in the proof of \textbf{Theorem~\ref{thm:weighted-objective}}
with $J=R$, these $4b^2$ choices are counted by
$\sum_{F\in\mathcal B}h(F,R)$. Thus
\[
 2mD^2\sum_{F\in\mathcal B}Q(F,R)=4b^2.
\]
Let $N_+$ count the choices for which the images of positions $1,2$
are adjacent in $G_{\mathrm s}$. These are the choices whose pattern
belongs to $\mathcal B_+$, so
$N_+=2mD^2\sum_{F\in\mathcal B_+}Q(F,R)$.

To count ordered pairs with a common endpoint, let $j_v$ be the
number of remaining marked edges incident with $v$. Choose a vertex
$v$ and an ordered pair $(\varepsilon_1,\varepsilon_2)$ of
remaining marked edges incident with $v$. Define
$\psi\colon P_3\to V(G_{\mathrm s})$ by $\psi(0)=v$ and set
$\psi(1)$ and $\psi(2)$ to the other endpoints of $\varepsilon_1$ and
$\varepsilon_2$, respectively. Take $\eta_{\{0,1\}}=\varepsilon_1$
and $\eta_{\{0,2\}}=\varepsilon_2$. The resulting tuple is counted
by $h(H_c,J_\star)$ for $c=c_G(\psi(1),\psi(2))$. Conversely, each
tuple counted by $h(H_c,J_\star)$ determines the choice
$(v,\varepsilon_1,\varepsilon_2)=(\psi(0),\eta_{\{0,1\}},
\eta_{\{0,2\}})$. Thus
\[
 2mD\sum_{c=0}^2 Q(H_c,J_\star)=\sum_v j_v^2.
\]
Let $a$ be the number of ordered pairs $(\varepsilon_1,\varepsilon_2)$
that share an endpoint. Each such pair is counted at most twice
in $\sum_v j_v^2$, once per common endpoint. Therefore
$a\geq\frac{1}{2}\sum_v j_v^2$.

For disjoint $\varepsilon_1,\varepsilon_2$, the strong clique
property gives adjacent vertices $u$ on $\varepsilon_1$ and $v$
on $\varepsilon_2$. Taking $(\phi(1),\phi(2))=(u,v)$ gives an
endpoint order counted by $N_+$. For a pair with common endpoint
$v$, write $x,y$ for the other endpoints of $\varepsilon_1$ and
$\varepsilon_2$, respectively. The choices
$(\phi(1),\phi(2))=(x,v)$ and $(v,y)$ are distinct and both are
counted by $N_+$, including when $x=y$. Thus the $b^2-a$ disjoint
pairs contribute at least one choice each, and the $a$ pairs
with a common endpoint contribute at least two each.
\textbf{Figure~\ref{fig:remaining-edge-orderings}} shows both cases.
It follows that
\[
 2mD^2\sum_{F\in\mathcal B_+}Q(F,R)
  =N_+\geq b^2+a\geq b^2+\frac{1}{2}\sum_v j_v^2.
\]
Multiplying by $\frac{2}{mD^2}$ and substituting the expressions
for $4b^2$ and $\sum_v j_v^2$ gives
\eqref{eq:remaining-edge-orderings}.
\end{proof}

The next example shows how a lower bound on $M$ leads to structural
restrictions and weighted density constraints. We use the resulting
counts to rule out $M\in\bigl[\frac{21}{16},\frac{11}{8}\bigr]$
when $D=4$.

\begin{theorem}\label{thm:degree-four-marked-edge}
Suppose $D=4$ and $M\geq\frac98$. If $u,v$ are joined by parallel
marked edges, then exactly two marked edges join them. Each of
$u,v$ has two neighbors outside $\{u,v\}$, the two neighbor sets
are disjoint, and all four edges to those neighbors are marked.
\end{theorem}

\begin{proof}
Since $D=4$, the lower bound on $M$ gives $m=16M\geq18$.
Suppose $k\geq2$ marked edges have endpoints $u,v$. Let
$A=(N_{G_{\mathrm s}}(u)\cup N_{G_{\mathrm s}}(v))\setminus\{u,v\}$
and put $q=|A|$. Each of $u,v$ is incident with at most $4-k$
edges to $A$, giving $q\leq2(4-k)$. At most $8-k$ marked
edges are incident with $u$ or $v$, since the sum of their marked
degrees is at most eight and counts the $k$ edges between them
twice.

Every marked edge disjoint from $\{u,v\}$ has an endpoint in $A$,
since it must be joined by an edge of $G$ to a marked edge with
endpoints $u,v$. Each vertex of $A$ has an edge to $u$ or $v$,
so it is incident with at most three edges disjoint from
$\{u,v\}$. Counting first the marked edges incident with $u$
or $v$ and then the others gives
\[
 18\leq m\leq8-k+3q\leq32-7k\leq18.
\]
All the inequalities are equalities. In particular, $k=2$, $q=4$,
and exactly six marked edges are incident with $u$ or $v$. After
the two marked edges between $u,v$ are counted, four marked edges
join $\{u,v\}$ to $A$. Each of $u,v$ has at most two edges to $A$.
Since $|A|=4$, each vertex of $A$ is incident with exactly one
of these four marked edges. Thus $u,v$ each have two neighbors
in $A$, their neighbor sets are disjoint, and all four edges are marked.
\end{proof}

\begin{theorem}\label{thm:degree-four-weighted-densities}
Suppose $D=4$ and $M\geq\frac98$. Let $F_2$ be the pattern on
$P_2$ with its pair of color $2$. Let $F_3$ be the pattern on
$P_3$ with $\{0,1\}$ and $\{1,2\}$ of color $2$ and $\{0,2\}$
of color $0$. Put $J=\{\{0,1\}\}$. Then
\begin{equation}\label{eq:degree-four-weighted-densities}
 4Q(F_3,J)=3Q(F_2,J).
\end{equation}
\end{theorem}

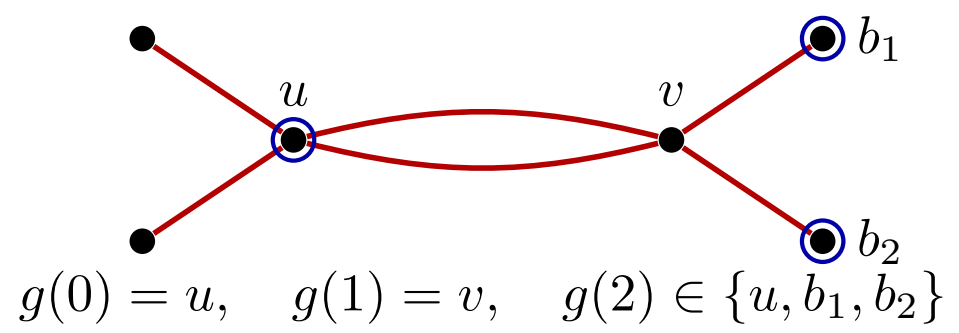
\begin{figure}[htbp]
\centering
\begin{tikzpicture}[every node/.style={font=\small},
                    marked edge/.append style={draw=red!70!black}]
  \node[graph vertex,label=above:$u$] (u) at (0,0) {};
  \node[graph vertex,label=above:$v$] (v) at (2.5,0) {};
  \node[graph vertex] (a1) at (-1,0.67) {};
  \node[graph vertex] (a2) at (-1,-0.67) {};
  \node[graph vertex,label=right:$b_1$] (b1) at (3.5,0.67) {};
  \node[graph vertex,label=right:$b_2$] (b2) at (3.5,-0.67) {};
  \draw[marked edge] (u) to[bend left=14] (v);
  \draw[marked edge] (u) to[bend right=14] (v);
  \draw[marked edge] (a1)--(u)--(a2);
  \draw[marked edge] (b1)--(v)--(b2);
  \foreach \point in {u,b1,b2}
    \draw[blue!65!black,thick] (\point) circle (3.9pt);
  \node at (1.25,-1.04) {$g(0)=u,\quad g(1)=v,\quad g(2)\in\{u,b_1,b_2\}$};
\end{tikzpicture}
\caption{The circled vertices are the three possible images of
position $2$ in \eqref{eq:degree-four-weighted-densities}.}
\label{fig:degree-four-marked-edge}
\end{figure}

\begin{proof}
Fix a tuple $(\phi,\eta_{\{0,1\}})$ counted by $h(F_2,J)$,
and write $u=\phi(0)$ and $v=\phi(1)$. The remaining marked
edge $\eta_{\{0,1\}}$ and the representative on $uv$ are
distinct. By \textbf{Theorem~\ref{thm:degree-four-marked-edge}},
$v$ has two neighbors $b_1,b_2$ other than $u$. Both edges
$vb_1,vb_2$ are red in $G_{\mathrm s}$, and neither $b_1$ nor
$b_2$ is adjacent to $u$.

To count extensions of this tuple, take
$g\colon P_3\to V(G_{\mathrm s})$ with $g|_{P_2}=\phi$ and keep
the edge $\eta_{\{0,1\}}$. A vertex $z$ gives a tuple counted by
$h(F_3,J)$ with $g(2)=z$ exactly when $c_G(v,z)=2$ and $c_G(u,z)=0$.
The first condition gives $z\in\{u,b_1,b_2\}$, as shown in
\textbf{Figure~\ref{fig:degree-four-marked-edge}}. All three
vertices satisfy the second condition, since $c_G(u,u)=0$
and $u$ is nonadjacent to $b_1,b_2$. Each tuple counted by
$h(F_3,J)$ has a unique restriction
$(g|_{P_2},\eta_{\{0,1\}})$ counted by $h(F_2,J)$, and each
tuple counted by $h(F_2,J)$ has three extensions. Hence
$h(F_3,J)=3h(F_2,J)$. Dividing by $2mD$ gives
$Q(F_3,J)=\frac{3}{D}Q(F_2,J)$, which is
\eqref{eq:degree-four-weighted-densities} because $D=4$.
\end{proof}

\begin{theorem}\label{thm:degree-four-interval}
If $D=4$, then $M\notin\bigl[\frac{21}{16},\frac{11}{8}\bigr]$.
\end{theorem}

\begin{proof}
Suppose $D=4$ and $M\in\bigl[\frac{21}{16},\frac{11}{8}\bigr]$.
We first prove $s\leq20$ to show that at least one of the
$m=16M\geq21$ marked edges is a remaining marked edge.
Let $H$ be the spanning subgraph
of $G_{\mathrm s}$ whose edges are the red edges, and write
$h_x=d_H(x)$ and $d_x=d_{G_{\mathrm s}}(x)$. Each $d_x$ is at
most $4$ because $G_{\mathrm s}$ records the distinct neighbors of
a vertex of $G$. For two disjoint red edges, their representatives
in $G$ belong to the marked strong clique. The endpoints of a joining
edge in $G$ are also adjacent in $G_{\mathrm s}$. Thus the edges
of $H$ form a strong clique in $G_{\mathrm s}$.

Suppose $s\geq21$. To show that $H$ has a vertex of degree $4$,
assume that all vertices have degree at most $3$ in $H$.
Fix a red edge $uv$ and put
$X=(N_{G_{\mathrm s}}(u)\cup N_{G_{\mathrm s}}(v))\setminus\{u,v\}$.
Every red edge other than $uv$ has an endpoint in $X$.
For an edge incident with $u$ or $v$, its other endpoint lies in
$X$. For an edge disjoint from $uv$, a joining edge to $uv$
ensures that one endpoint lies in $X$. Both $u$ and $v$ have at
most three neighbors in $X$, so $|X|\leq6$. It follows that
$s-1\leq\sum_{x\in X}h_x\leq18$, a contradiction.
Fix a vertex $v$ with $h_v=4$.

Suppose $H$ contains a triangle on vertex set $T$. For each vertex
$x$, put $a_x=|N_{G_{\mathrm s}}(x)\cap T|$. For every red edge
$xy$, the set $(N_{G_{\mathrm s}}(x)\cup N_{G_{\mathrm s}}(y))\cap T$
meets every triangle edge. A triangle edge disjoint from $xy$
has a joining edge to $xy$, and if the two edges share an endpoint,
the other endpoint of the triangle edge belongs to the set.
Thus $a_x+a_y\geq2$. If exactly one endpoint of $xy$ belongs to
$T$, it has two neighbors in $T$ and the other endpoint has at
least one. If both endpoints belong to $T$, each has two neighbors
there. Hence $a_x+a_y\geq2+|\{x,y\}\cap T|$. Summing over
$xy\in E(H)$ gives the first inequality below.
\[
 2s+\sum_{x\in T}h_x
 \leq\sum_x a_xh_x
 \leq4\sum_{x\in T}d_x-\sum_{x\in T}(d_x-h_x).
\]
For the second inequality, write
$\sum_x a_xh_x=\sum_x a_xd_x-\sum_x a_x(d_x-h_x)$.
The first sum is at most $4\sum_{x\in T}d_x$, since $d_x\leq4$
and $\sum_x a_x=\sum_{x\in T}d_x$. The second sum counts pairs
consisting of a neighbor of $x$ in $T$ and a blue edge incident
with $x$. For each blue incidence at a vertex $y\in T$ along an edge
$xy$, take the pair $(y,xy)$ at $x$. These choices are distinct, so
$\sum_x a_x(d_x-h_x)\geq\sum_{y\in T}(d_y-h_y)$.
Combining the displayed bounds gives $2s\leq3\sum_{x\in T}d_x\leq36$,
a contradiction. Hence $H$ is triangle-free.

Let $A=N_{G_{\mathrm s}}(v)$ and
$W=V(G_{\mathrm s})\setminus(A\cup\{v\})$. Let $t$ count the red
edges between $A$ and $W$, and let $w$ count the red edges within
$W$. The four edges at $v$ are red, so $|A|=4$. A red edge within
$A$ would give a red triangle with $v$. Thus $s=4+t+w$. For
$x\in W$, put $c_x=|N_{G_{\mathrm s}}(x)\cap A|$ and
$p=\sum_{x\in W}c_x$. This sum counts all edges between $A$ and
$W$, so $t\leq p$. Each vertex of $A$ is adjacent to $v$ and has
at most three neighbors in $W$, giving $p\leq12$.

For each red edge $xy$ within $W$ and each $a\in A$, the red edges
$xy$ and $va$ must have a joining edge. Neither $x$ nor $y$ is
adjacent to $v$, so $a$ is adjacent to at least one of them.
Thus $c_x+c_y\geq4$. The edge $xy$ uses one incidence at each
endpoint, so $c_x,c_y\leq3$ and both are positive. Summing
$c_x+c_y-2$ over the red edges within $W$ gives
\[
 2w\leq
 \sum_{\substack{x\in W\\c_x\geq1}}(c_x-1)d_{H[W]}(x)
 \leq p-\sum_{\substack{x\in W\\c_x\geq1}}(c_x-2)^2
 \leq12,
\]
where the second inequality follows from
$d_{H[W]}(x)\leq4-c_x$ and
$(c_x-1)(4-c_x)=c_x-(c_x-2)^2$.
If $w=6$, equality throughout gives $p=12$, $c_x=2$ whenever
$c_x>0$, and $d_{H[W]}(x)=4-c_x=2$ at each of these vertices.
There are six such vertices, and a triangle-free simple graph
with six vertices of degree $2$ and all other vertices isolated
has a six-cycle. Each cycle vertex has two neighbors on the
cycle and two in $A$, so there is no chord in $G_{\mathrm s}$.
Opposite edges of the cycle then have no joining edge in
$G_{\mathrm s}$, a contradiction.

We therefore have $w\leq5$. Since $s\geq21$, equality holds
in $s=4+t+w\leq4+12+5$, giving $t=p=12$ and $w=5$. All edges
between $A$ and $W$ are red. For every red edge within $W$,
its endpoints have complementary neighbor sets in $A$, since
their neighbor sets cover $A$ and a common neighbor would give
a red triangle.

There can be no red four-cycle in $W$ on vertices with $c_x=2$.
Each vertex on such a cycle would have two neighbors on the cycle
and two in $A$. The fifth red edge of $H[W]$ would therefore be
disjoint from the cycle and have no joining edge to any cycle edge.

Suppose $\ell\geq1$ of the five red edges have one endpoint with
$c_x=1$ and one with $c_x=3$. The $\ell$ endpoints with $c_x=3$
are distinct because each has at most one neighbor in $W$.
At least one vertex has $c_x=1$. Each of the other $5-\ell$ red
edges has both endpoints with $c_x=2$. Each such vertex has at
most two neighbors in $W$, so these edges require at least
$5-\ell$ vertices with $c_x=2$. Thus
\[
 12=p\geq3\ell+1+2(5-\ell),
\]
This forces $\ell=1$ and equality in each of the three vertex
counts. In particular, the four edges between vertices with
$c_x=2$ lie on exactly four vertices, each of degree at most
$2$ in $H[W]$. All four degrees must therefore equal $2$.
They form a red four-cycle, a contradiction.

It follows that all five red edges within $W$ join vertices with
$c_x=2$. There are at most six such vertices, each of degree at
most $2$ in $H[W]$. Along these red edges the neighbor sets in
$A$ alternate between complementary sets, so there is no odd
red cycle. We have already excluded a red four-cycle on these
vertices, and a longer even cycle requires at least six edges. The five
red edges therefore form a forest on at most six vertices
incident with them. It has six vertices and is connected, so
its degree bound of $2$ makes it a path $x_1,\ldots,x_6$.
Both $x_2$ and $x_3$ have two path neighbors and two neighbors
in $A$, leaving no edge joining $x_2x_3$ to $x_5x_6$.
This last contradiction gives $s\leq20$.

Since $m\geq21$ and $s\leq20$, we have $b=m-s\geq1$ and $e>0$.
To use the condition that at least one remaining marked edge
exists, we count pairs of marked edges whose first edge is a
remaining marked edge, with an endpoint order for each edge. Put
\[
 L=\sum_{F\in\mathcal B}\sum_{\substack{J\subseteq R\\ \{0,1\}\in J}}Q(F,J).
\]
In the counting correspondence of
\textbf{Theorem~\ref{thm:weighted-objective}}, the condition
$\{0,1\}\in J$ selects precisely these pairs. By
\eqref{eq:weighted-raw-count}, their number is $2mD^2L$.
There are $2b$ choices for the first edge with its
endpoint order and $2m$ for the second. The two choices are
independent and may use the same edge.
Consequently
\[
 L=\frac{(2b)(2m)}{2mD^2}=2Me\geq\frac{21}{8}e.
\]

To bound $L$ from above, fix the first edge and its endpoint order.
This edge and the representative with the same endpoints are
distinct marked edges. Thus the count in the proof of
\textbf{Theorem~\ref{thm:degree-four-marked-edge}} gives $m\leq18$.
Thus there are at most $36$ choices of the second marked edge
with its endpoint order. It follows that
\[
 L\leq\frac{(2b)36}{2mD^2}=\frac{9}{4}e.
\]
Since $\frac{21}{8}>\frac94$ and $e>0$, the two bounds on $L$
contradict each other.
\end{proof}

For the other intervals of $M$, we give the arguments ruling them out
as part of the complete proof of \textbf{Conjecture~\ref{conj:fgst}}
in \textbf{Sections~\ref{sec:uniform} and~\ref{sec:exact}}.
The derivations use the same counting ideas as the examples above.
The larger systems of weighted density constraints require extensive
case analysis and exact computational verification.
Hence the proof is computer-assisted. In
\textbf{Section~\ref{sec:uniform}}, we work through one interval,
from the structural analysis to a weighted density constraint and
then to an upper bound on $M$ that excludes the interval. For the
other intervals and degree ranges,
\textbf{Sections~\ref{sec:uniform} and~\ref{sec:exact}} state the
constraints and estimates and refer to the proofs in
\textbf{Sections~S3 and~S4} and the checks in
\textbf{Appendices~SG.3 and~SG.4} of~\cite{FGSTSupplement}.
The project repository~\cite{FGSTCode} contains the companion
manuscript in \path{manuscript/} and the verification programs and
data in \path{reproducibility/}.

\section{The Proof of the Bound \texorpdfstring{$\frac{5}{4}\Delta^2$}{5 Delta squared / 4}}
\label{sec:uniform}

In this section, we prove \textbf{Theorem~\ref{prop:multi-uniform}}, which
gives a uniform upper bound on the size of a strong clique in a loopless
multigraph. We first narrow the range $M>\frac{5}{4}$, then partition
the remaining range and rule out each interval.

\begin{theorem}\label{prop:multi-uniform}
Let $D\geq0$ be an integer and let $G$ be a finite loopless multigraph
with $\Delta(G)\leq D$. Every strong clique in $G$ has at most
$\frac{5}{4}D^2$ edges.
\end{theorem}

\subsection{Partition of the Range of \texorpdfstring{$M$}{M}}
\label{subsec:uniform-partition}

Fix a strong clique of $m>0$ edges, so $D>0$. Mark its edges and
reduce the multigraph as in \textbf{Subsection~\ref{subsec:weighted-sketch}}.
Write $G_{\mathrm s}$ for the underlying simple graph of the
reduced multigraph and $s$ for its number of red edges. The red edges
form a strong clique in $G_{\mathrm s}$, whose maximum degree is at
most $D$. Thus $S=\frac{s}{D^2}\leq\frac{5}{4}$ by
\textbf{Lemma~S3.7} of~\cite{FGSTSupplement}.
When taking weighted densities, use the $D$-regular multigraph
obtained by applying \textbf{Lemma~\ref{lem:regularization}} to
the reduced multigraph. Its new edges are all unmarked, and its
underlying simple graph has the same $s$ red edges.
In both graphs $e=\frac{m-s}{m}$, so $S=M(1-e)$.

To restrict the possible values of $M$, suppose $M>\frac{5}{4}$.
Five vertices with all ten pairs red
would give $4m\leq5D(D-1)-20$ by
(S3.27) of~\cite{FGSTSupplement}, contrary to
$4m>5D^2$. Thus every pattern containing a red $K_5$ has density
zero. The five-position checks in \textbf{Appendix~SG.4}
of~\cite{FGSTSupplement} then give $M<1.28746$ and
$e<\frac{421}{2000}$. To improve the upper bound on $M$, the
same appendix applies the certificates
\path{baseline_half_refit.json}, \path{direct_half_refit.json}
and \path{direct_half_scs.json} in succession. The rational
comparisons there give (SG.1), namely
\begin{equation}\label{eq:uniform-initial-range}
 M<1.28492,\qquad e<\frac{421}{2000}.
\end{equation}
We partition this range of $M$ as follows.
\begin{equation}\label{eq:uniform-partition}
\begin{aligned}
 I_1&=[1.283,1.28492), & I_2&=[1.279,1.283),\\
 I_3&=[1.274,1.279), & I_4&=[1.268,1.274),\\
 I_5&=\left(\frac{5}{4},1.268\right). &&
\end{aligned}
\end{equation}
For each interval, let $M_-$ and $M_+$ be its left and right
endpoints. We impose $M_-\leq M\leq M_+$ when deriving bounds for
that interval. An upper bound strictly below $M_-$ then excludes
the interval. All terminating decimals in this section denote exact
rational numbers.

\subsection{Structural Analysis in Each Interval}\label{subsec:uniform-structure}

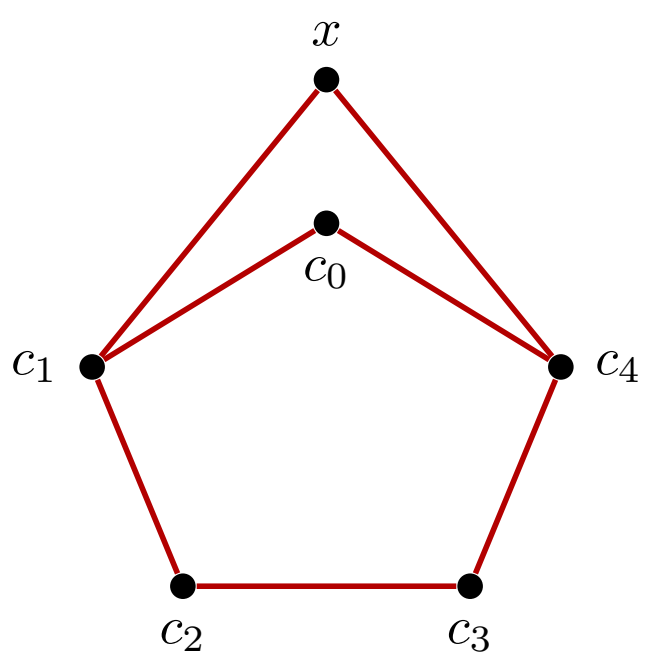
\begin{figure}[htbp]
\centering
\begin{tikzpicture}[every node/.style={font=\small},
                    marked edge/.append style={draw=red!70!black}]
  \node[graph vertex,label=below:$c_0$] (c0) at (0,1.4) {};
  \node[graph vertex,label=left:$c_1$] (c1) at (-1.55,0.45) {};
  \node[graph vertex,label=below:$c_2$] (c2) at (-0.95,-1) {};
  \node[graph vertex,label=below:$c_3$] (c3) at (0.95,-1) {};
  \node[graph vertex,label=right:$c_4$] (c4) at (1.55,0.45) {};
  \node[graph vertex,label=above:$x$] (x) at (0,2.35) {};
  \draw[marked edge] (c0)--(c1)--(c2)--(c3)--(c4)--(c0);
  \draw[marked edge] (c1)--(x)--(c4);
\end{tikzpicture}
\caption{A vertex $x\in U_0$ has exactly two cycle neighbors, and
both edges from $x$ to the cycle are red.}
\label{fig:main-uniform-cycle-parts}
\end{figure}

For the structural counts, work in the reduced multigraph. An
induced red five-cycle is a cycle of length five in $G_{\mathrm s}$
whose cycle edges are red and which has no chord in $G_{\mathrm s}$.
Fix an induced red five-cycle $C=c_0c_1c_2c_3c_4c_0$, with indices
modulo five. For each $i\in P_5$, let $U_i$ consist of the vertices
$x$ with $N_{G_{\mathrm s}}(x)\cap V(C)=\{c_{i-1},c_{i+1}\}$
and with both edges $xc_{i-1},xc_{i+1}$ red.
These five sets are disjoint, and $c_i\in U_i$. For example,
\textbf{Figure~\ref{fig:main-uniform-cycle-parts}} shows a vertex of $U_0$.
Put $B_C=\bigcup_i U_i$. To bound the number of remaining incidences
at the cycle vertices, count one edge to each neighbor of $c_i$
in $B_C$. These neighbors are exactly $U_{i-1}\cup U_{i+1}$, so
at most $D-|U_{i-1}|-|U_{i+1}|$ further edges are incident with
$c_i$. Summing these five bounds gives
\begin{equation}\label{eq:main-uniform-cycle-capacity}
 R(C)=\sum_{i=0}^4\bigl(D-|U_{i-1}|-|U_{i+1}|\bigr)=5D-2|B_C|.
\end{equation}

Every vertex of $U_i$ is adjacent to every vertex of $U_{i+1}$.
Indeed, for $x\in U_i$ and $y\in U_{i+1}$, the red edges
$xc_{i-1}$ and $yc_{i+2}$ are disjoint. The definitions of the
sets exclude $xc_{i+2}$ and $yc_{i-1}$, and the cycle has no
chord $c_{i-1}c_{i+2}$. Hence $xy$ must be an edge. Let $L$ be
the number of blue edges between consecutive sets. For each
vertex $x$, write $h_x$ for the number of marked edges incident
with $x$. For $x\notin B_C$, let $b_x$ count the blue edges
from $x$ to $B_C$. The edges counted in $L$ use $2L$ unmarked
incidences in $B_C$, and the blue edges from outside $B_C$ use
another $\sum_{x\notin B_C}b_x$ unmarked incidences there.
Thus the sum of marked degrees in $B_C$ is at most
$D|B_C|-2L-\sum_{x\notin B_C}b_x$. Adding the marked degrees
outside $B_C$ gives
\begin{equation}\label{eq:main-uniform-core-cut}
 2m\leq D|B_C|-2L+\sum_{x\notin B_C}(h_x-b_x).
\end{equation}
To count the adjacencies to the cycle, define the
\emph{signature} of every vertex $x\in V(G)$ to be
$\operatorname{sig}(x)=\{i\in P_5:xc_i\in E(G_{\mathrm s})\}$.
Put $K=\sum_{x\notin B_C}|\operatorname{sig}(x)|$. For each
adjacency counted in $K$, there is an edge of $G$ from
$x\notin B_C$ to one of the cycle vertices. These incidences
are included in the bound $R(C)$, so $K\leq R(C)$.

For each marked edge $xy$, the signature union
$\operatorname{sig}(x)\cup\operatorname{sig}(y)$ contains an
endpoint index of every edge of $C$. If $xy$ and a cycle edge
are disjoint, the strong clique condition gives an adjacency
from $x$ or $y$ to an endpoint of the cycle edge. If they share
an endpoint, that vertex is adjacent to the other endpoint of
the cycle edge. Thus the cycle vertices with indices in the union
form a vertex cover of $C$, so the union has size at least three.
In particular, a marked
edge from a vertex with empty signature must have its other
endpoint at a vertex with at least three cycle neighbors.

Put $r=\frac{R(C)}{D}$ and $g=M-\frac54$. For $R(C)\leq D$,
the signature counts in (S3.24) of~\cite{FGSTSupplement} give
$\sum_{x\notin B_C}(h_x-b_x)\leq\max\{\frac D2,R(C)\}K$.
If $R(C)\leq\frac D2$, substitute this bound and $K\leq R(C)$
into \eqref{eq:main-uniform-core-cut} to obtain
$2m\leq D|B_C|+\frac D2R(C)=\frac52D^2$ by
\eqref{eq:main-uniform-cycle-capacity}. This contradicts
$M>\frac54$, so $r>\frac12$.
For $\frac12\leq r\leq1$, the same estimate and $K\leq R(C)$ give
$2m\leq D|B_C|+R(C)^2$. Dividing by $D^2$ and substituting
$|B_C|=\frac{5D-R(C)}{2}$ gives the first inequality below.
The second follows from the refined
count (S3.22) of~\cite{FGSTSupplement} by dropping
the terms contributed by the blue edges between consecutive sets
and the remaining marked edges at the cycle vertices.
\begin{equation}\label{eq:main-uniform-r-bounds}
 2g\leq r\left(r-\frac12\right),\qquad
 2g\leq\frac{1+r}{2+r}\,r\left(r-\frac12\right).
\end{equation}
Both right-hand sides in \eqref{eq:main-uniform-r-bounds}
are strictly increasing on $[\frac12,1]$. In an interval with
left endpoint $M_-$, we have $2g\geq2(M_--\frac54)$. Therefore
if one of the right-hand sides at $\alpha\in[\frac12,1]$ is
at most $2(M_--\frac54)$, then $r\geq\alpha$.
For $I_1$, the first estimate gives
\[
 0.6084\left(0.6084-\frac12\right)=0.06595056
 <2\left(1.283-\frac54\right)=0.066.
\]
Consequently $r>0.6084$ in $I_1$. Applying the second estimate
at the left endpoints of $I_2$, $I_3$ and $I_4$ gives
$r\geq0.6446$, $r\geq0.6242$ and $r\geq0.5978$, respectively.
In $I_5$ we use $r>\frac12$. Each stated lower bound is at most
one, so it also holds when $r>1$.

\subsection{Weighted Density Constraints in Each Interval}
\label{subsec:uniform-weighted-constraints}

Let $\mathcal C$ be the set of ordered tuples $(c_0,\ldots,c_4)$
that define induced red five-cycles $c_0c_1c_2c_3c_4c_0$ in
$G_{\mathrm s}$. Each cycle has five choices for $c_0$ and
two choices for the direction, so it gives ten tuples.
Define the normalized count $q_C$ and the normalized sum $L_R$ by
\begin{equation}\label{eq:ordered-cycle-densities}
 q_C=\frac{|\mathcal C|}{2sD^3},\qquad
 L_R=\frac{\sum_{C\in\mathcal C}\frac{R(C)}{D}}{2sD^3}.
\end{equation}
Whenever $\frac{R(C)}{D}\geq\alpha$ for every $C\in\mathcal C$,
summing this inequality gives $L_R\geq\alpha q_C$. Both sides
are zero if $\mathcal C$ is empty. For each interval in
\eqref{eq:uniform-partition}, we use the lower bound $\alpha$
from \textbf{Subsection~\ref{subsec:uniform-structure}}.

To write this constraint in weighted densities, let $F_C$ be
the pattern on $P_5$ with red pairs $\{i,i+1\}$ for $i\in P_5$,
with indices modulo five, and color $0$ on the other pairs.
A map giving the colors of $F_C$ is
injective. Positions on a red pair have distinct images. For
nonadjacent positions $i,j$, there is a position $k$ for which
one of $\{i,k\},\{j,k\}$ is red and the other has color $0$.
If the images of $i,j$ coincide, those two pairs must receive
the same color, a contradiction. Thus the map sending $i$ to
$c_i$ gives a bijection between the maps and $\mathcal C$.

Let $\mathcal F_B$ consist of the extensions of $F_C$ to $P_6$
in which the pairs from position $5$ to $i-1$ and $i+1$ are
red for some $i\in P_5$, with indices modulo five, and all
other pairs from position $5$ to the cycle have color $0$.
For a fixed cycle map, the extensions to these patterns are
exactly the choices of the image of position $5$ in $B_C$.
This choice may be a vertex of the cycle, since the maps in
\eqref{eq:weighted-density} may repeat images.
All images in both counts are incident with red edges, so
they belong to the reduced multigraph and their pair colors
are preserved by regular completion. With $J=\varnothing$,
the five- and six-position denominators in
\eqref{eq:weighted-density} are $2mD^3$ and $2mD^4$,
respectively. Hence
\[
 Q(F_C,\varnothing)=\frac{|\mathcal C|}{2mD^3},\qquad
 \sum_{F\in\mathcal F_B}Q(F,\varnothing)
    =\frac{\sum_{C\in\mathcal C}|B_C|}{2mD^4}.
\]
Substituting $R(C)=5D-2|B_C|$ and $\frac{s}{m}=1-e$ into
\eqref{eq:ordered-cycle-densities} gives
\begin{equation}\label{eq:main-uniform-cycle-constraint}
 (1-e)(L_R-\alpha q_C)
  =(5-\alpha)Q(F_C,\varnothing)
      -2\sum_{F\in\mathcal F_B}Q(F,\varnothing)\geq0.
\end{equation}
We use the simple-graph certificates to bound $S$ and hence
to bound $e=1-\frac{S}{M}$ from below. Complete $G_{\mathrm s}$ to a
$D$-regular simple graph by \textbf{Lemma~\ref{lem:regularization}}
and color the new edges blue. Every vertex of a red cycle or
of one of its sets $U_i$ is incident with a red edge, so this
completion preserves the counts in
\eqref{eq:ordered-cycle-densities}. In the simple-graph
certificates the normalization uses $s$ marked edges, which
is why the inequalities below use $q_C$ and $L_R$ as defined
in \eqref{eq:ordered-cycle-densities}.
Equation~(S3.26) of~\cite{FGSTSupplement}
gives
\begin{equation}\label{eq:main-uniform-simple-certificates}
\begin{aligned}
 S-0.99999q_C
   &\leq \frac54-\frac{0.99999}{8}+2\cdot10^{-9},\\
 S+0.5682L_R
   &\leq \frac54+10^{-9},\\
 S+0.4897(L_R-0.32q_C)
   &\leq \frac54-\frac{0.4897\cdot0.32}{8}+2\cdot10^{-9},\\
  S+0.45119\left(L_R-\frac{q_C}{2}\right)
   &\leq \frac54-\frac{0.45119}{16}+2\cdot10^{-9}.
\end{aligned}
\end{equation}
The certificate inputs and their checks are given in
\textbf{Appendix~SG.1} of~\cite{FGSTSupplement}.
Substitute $L_R\geq\alpha q_C$ into the last three inequalities.
Each of the four inequalities now bounds $S$ from above by an
affine function of $q_C$. Include the constant bound
$S\leq\frac54$ and let $S_{\max}(\alpha)$ be the maximum, over
$q_C\geq0$, of the minimum of these five functions. This minimum
is piecewise affine and tends to $-\infty$ because one of the
functions has slope $-0.5682\alpha<0$. Its maximum is attained
at zero or at an intersection of two of the lines. Evaluating
at those points uses rational arithmetic. Since $M\geq M_->0$
and $S\leq S_{\max}(\alpha)$, we obtain
\[
 e=1-\frac{S}{M}\geq 1-\frac{S_{\max}(\alpha)}{M_-}.
\]
\textbf{Table~\ref{tab:main-uniform-intervals}} records the resulting
lower bounds on $e$.

For the final interval, $\alpha=\frac12$, so the last line of
\eqref{eq:main-uniform-simple-certificates} alone gives
\[
 S\leq\frac54-\frac{0.45119}{16}+2\cdot10^{-9}=1.221800627.
\]
Since $M>\frac{5}{4}$, this gives
$e\geq1-\frac{4\cdot1.221800627}{5}=0.0225594984>0.0225$.

\begin{lateproof}{\textbf{Theorem~\ref{prop:multi-uniform}}}
The size bound holds when $m=0$. Suppose that $m>0$ and
$M>\frac54$. Then $M$ belongs to one of the intervals in
\eqref{eq:uniform-partition}. \textbf{Subsection~\ref{subsec:uniform-weighted-constraints}}
gives the lower bound on $e$ in \textbf{Table~\ref{tab:main-uniform-intervals}},
and \eqref{eq:uniform-initial-range} gives $e<\frac{421}{2000}$.
We first use these bounds for $I_5$, where
$\frac{9}{400}\leq e<\frac{421}{2000}$ and $\frac54<M<1.268$.

Let $x$ have one coordinate for each isomorphism class of
weighted patterns used in this calculation. For a representative
$(F,J)$, the coordinate is $\frac{Q(F,J)}{|\operatorname{Aut}(F,J)|}$,
where $\operatorname{Aut}(F,J)$ is the group of permutations of
the positions preserving all pair colors and the set $J$.
Thus $x_j\geq0$ for every coordinate. Isomorphic patterns
have equal weighted densities, since relabeling the positions
gives a bijection of their maps and choices of remaining
marked edges. Consequently \eqref{eq:weighted-objective}
expresses $M$ as a linear form in $x$.

The verifier writes the counting equations as $Ax=a$, where
$a=(1,0,\ldots,0)^{\mathsf T}$. The equation with right-hand
side $1$ is the normalization in \eqref{eq:weighted-normalization}.
In the linear inequalities, replace each constant by that
constant times the left-hand side of this normalization. For
example, the upper bound $M\leq1.274$ can be written as
\[
 1.274\bigl(Q(F_\bullet,\varnothing)+Q(F_\bullet,\{p_\bullet\})\bigr)
 -M\geq0.
\]
In this form the inequalities are linear forms $g_\nu(x)\geq0$.
The system for $I_5$ uses $\frac54\leq M\leq1.274$ and
$\frac{9}{400}\leq e\leq\frac{421}{2000}$, which contain the
required ranges, and includes
\eqref{eq:main-uniform-cycle-constraint}.
Write $W_i(x)$ for its positive semidefinite matrices,
whose entries are linear in $x$. These include the weighted
flag matrices and the matrices from the underlying simple
graph expressed with the weighted normalization.
\textbf{Appendices~SG.1, SG.3 and~SG.4}
of~\cite{FGSTSupplement} give the additional constraints and the
indexed system.

The file \path{rank6_uniform_sharp_125_exclusion.json} in
the repository~\cite{FGSTCode} records the rational multipliers
for the equations, inequalities and matrix constraints. Use
$\lambda$ for the vector of equation multipliers. For the
inequalities, it gives nonnegative multipliers $\rho_\nu$
for $g_\nu(x)$. For the matrix constraints, it records integer
vectors $z_{ia}$. The associated quadratic forms
$z_{ia}^{\mathsf T}W_i(x)z_{ia}$ have weight
$w_{ia}=10^{-18}$, the reciprocal square of the stored
vector denominator.
Define the residual coefficients $r_j$ by the identity of
linear forms in $x$
\[
 \sum_j r_jx_j=\lambda^{\mathsf T}Ax-M
   -\sum_\nu\rho_\nu g_\nu(x)
   -\sum_{i,a}w_{ia}\,z_{ia}^{\mathsf T}W_i(x)z_{ia}.
\]
For any $x$ satisfying the counting equations $Ax=a$, we have
$\lambda^{\mathsf T}Ax=\lambda^{\mathsf T}a=C$, where $C$ is
the multiplier of the normalization equation.

To bound $\sum_{j:r_j<0}(-r_j)x_j$, group the coordinates
by the number of positions and the size of $J$. For a group
$I$, put $a_j=\max\{0,-r_j\}$ for $j\in I$. We bound
$\sum_{j\in I}a_jx_j$ in two ways. First use the individual
upper bounds on $x_j$ from \eqref{eq:weighted-density-bounds},
divided by the automorphism group sizes. Second use the bound
on $\sum_{j\in I}x_j$ multiplied by $\max_{j\in I}a_j$.
Both upper bounds are constant multiples of $1-e$ when $J$
is empty and of $e$ otherwise. Taking the smaller bound in
each group and summing therefore gives an affine function of $e$.
Its maximum on
$[\frac{9}{400},\frac{421}{2000}]$ is attained at an endpoint.
Call this maximum $B$. This calculation is given with
(SG.2) of~\cite{FGSTSupplement}.
Substituting $C$ and rearranging the coefficient identity gives
\begin{equation}\label{eq:main-uniform-final-decomposition}
 C+B-M
 =\sum_\nu \rho_\nu g_\nu(x)
   +\sum_{i,a}w_{ia}\,z_{ia}^{\mathsf T}W_i(x)z_{ia}
   +\left(B+\sum_j r_jx_j\right).
\end{equation}
The first sum is nonnegative because $\rho_\nu\geq0$ and
$g_\nu(x)\geq0$. The second is nonnegative because
$w_{ia}\geq0$ and $W_i(x)\succeq0$. The last term is at least
$B-\sum_{j:r_j<0}(-r_j)x_j\geq0$ because every $x_j\geq0$.
The exact values from (SG.2)
and its derivation in~\cite{FGSTSupplement} are
\[
 C=\frac{312419970929}{250000000000},
 \qquad B=\frac{73834315704703985647}{100000000000000000000000000}.
\]
Thus \eqref{eq:main-uniform-final-decomposition} gives
\begin{equation}\label{eq:main-uniform-final-certificate}
 M\leq C+B=\frac{124968062205915704703985647}
                   {100000000000000000000000000}
       <\frac{5}{4}.
\end{equation}
The bound lies below the left endpoint of $I_5$, as shown in
\textbf{Figure~\ref{fig:uniform-final-interval}}.

\begin{figure}[htbp]
\centering
\begin{tikzpicture}[font=\small]
  \draw[gray,->] (0,0)--(9.5,0) node[right] {$M$};
  \draw[red!65!black,line width=1.1pt,<-] (0.4,0)--(2.7,0);
  \draw[blue!65!black,line width=1.1pt] (4.65,0)--(8.35,0);
  \draw (2.7,-0.09)--(2.7,0.09);
  \node[below] at (2.7,-0.11) {$C+B$};
  \filldraw[draw=blue!65!black,fill=white] (4.65,0) circle (2pt);
  \filldraw[draw=blue!65!black,fill=white] (8.35,0) circle (2pt);
  \node[below] at (4.65,-0.11) {$\frac{5}{4}$};
  \node[below] at (8.35,-0.11) {$1.268$};
  \node[above] at (1.5,0.12) {$M\leq C+B$};
  \node[above] at (6.5,0.12) {$I_5$};
\end{tikzpicture}
\caption{The upper bound from \eqref{eq:main-uniform-final-decomposition}
lies below the left endpoint of the assumed interval $I_5$.}
\label{fig:uniform-final-interval}
\end{figure}
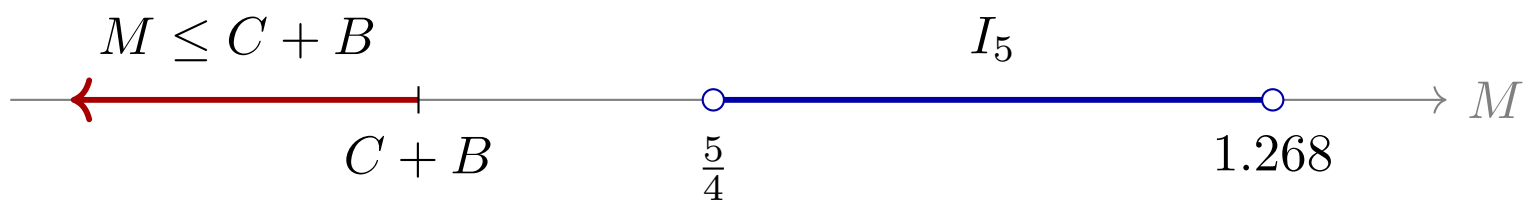

The verifier \path{verify_shadow_chain.py} repeats this
coefficient comparison with the multipliers for $I_1,\ldots,I_4$,
the first four intervals in \textbf{Table~S3.1} of~\cite{FGSTSupplement}.
\textbf{Appendix~SG.4} of~\cite{FGSTSupplement}
lists the exact rational bounds. The displayed bounds in
\textbf{Table~\ref{tab:main-uniform-intervals}} are rounded upwards and
are all below their respective lower endpoints $M_-$. This
excludes every interval in \eqref{eq:uniform-partition}, proving
$m\leq\frac54D^2$.
\end{lateproof}

\begin{table}[htbp]
\centering
\begin{tabular}{@{}lrrr@{}}
\toprule
Interval in \eqref{eq:uniform-partition} & $\alpha$ & Lower bound on $e$ & Upper bound on $M$\\
\midrule
$I_1$ & $0.6084$ & $0.051$ & $1.274527$\\
$I_2$ & $0.6446$ & $0.049$ & $1.269533$\\
$I_3$ & $0.6242$ & $0.0445$ & $1.262888$\\
$I_4$ & $0.5978$ & $0.0395$ & $1.264842$\\
$I_5$ & $\frac{1}{2}$ & $0.0225$ & $1.249681$\\
\bottomrule
\end{tabular}
\caption{The lower bounds $\frac{R(C)}{D}\geq\alpha$, lower bounds
on $e$, and rounded upper bounds on $M$ for the intervals in
\eqref{eq:uniform-partition}.}
\label{tab:main-uniform-intervals}
\end{table}

For even $D$, the bound in \textbf{Theorem~\ref{prop:multi-uniform}}
equals $f(D)$.

\section{The Proof of the Bound \texorpdfstring{$\Delta^2+\lfloor\frac{\Delta}{2}\rfloor^2$}{Delta squared plus floor(Delta/2) squared}}
\label{sec:exact}\label{sec:multi-odd}

In this section, we complete the proof of \textbf{Conjecture~\ref{conj:fgst}}.
Let $G$ be a finite loopless multigraph and put $D=\Delta(G)$. If $D$ is
even, then $f(D)=\frac{5}{4}D^2$, so
\textbf{Theorem~\ref{prop:multi-uniform}} gives the required bound.
For odd $D$, we use part \textup{(iii)} of
\textbf{Theorem~\ref{thm:weighted-objective}} to reduce the proof to
ruling out a strong clique with exactly $f(D)+1$ edges.

\subsection{The Structure of a Counterexample}\label{subsec:exact-structure}

Suppose $D\geq3$ is odd and a strong clique in $G$ has more than
$f(D)$ edges. Mark exactly $m=f(D)+1$ of these edges and apply the
reduction in \textbf{Subsection~\ref{subsec:weighted-sketch}}.
Write $G$ for the reduced multigraph and $G_{\mathrm s}$ for its
underlying simple graph. The $s$ red edges of $G_{\mathrm s}$
form a strong clique. There are $b=m-s$ remaining marked edges.
Since $f(D)=\frac{5D^2-2D+1}{4}$ for odd $D$, the quantities
$M=\frac{m}{D^2}$, $S=\frac{s}{D^2}$ and $e=\frac{b}{m}$ satisfy
\begin{equation}\label{eq:odd-degree-mass}
 M=\frac{5}{4}-\frac{1}{2D}+\frac{5}{4D^2},
 \qquad S=M(1-e).
\end{equation}

To show that $G_{\mathrm s}$ contains an induced red five-cycle
when $D\geq9$, subtract $\frac{98}{81}$ from
\eqref{eq:odd-degree-mass} to obtain
\[
 M-\frac{98}{81}=\frac{(D-9)(13D-45)}{324D^2}\geq0.
\]
Equation (S3.28) of~\cite{FGSTSupplement} gives
$M<\frac{98}{81}$ if there is no induced red five-cycle. Thus
$G_{\mathrm s}$ contains one. For each induced red five-cycle $C$,
let $U_i$ and $R(C)$ be as in
\textbf{Subsection~\ref{subsec:uniform-structure}}. Since
$R(C)=5D-2|\bigcup_i U_i|$ and $D$ is odd, $R(C)$ is odd.
\textbf{Lemma~S4.1} of~\cite{FGSTSupplement} gives
$R(C)>\frac{D+1}{2}$ for every induced red five-cycle $C$.

We next show that $b>0$, so two marked edges have the same
endpoints. For $D=3,5,7,9$, the exact simple-graph bound in
\textbf{Lemma~S4.8} of~\cite{FGSTSupplement}
applied to the red edges gives $s\leq f(D)<m$.
For $D\geq11$, let $q_C,L_R$ be as in
\eqref{eq:ordered-cycle-densities}. The bound on $R(C)$ gives
$L_R\geq\frac{q_C}{2}$, so the last inequality in
\eqref{eq:main-uniform-simple-certificates} gives
$S\leq1.221800627$. For $D\geq15$,
\[
 M-\frac{11}{9}=\frac{(D-15)(D-3)}{36D^2}\geq0,
\]
and hence $S<M$.

For $D=11,13$, \textbf{Lemma~S4.2}
of~\cite{FGSTSupplement} gives $R(C)\geq9$ in both cases.
Put $\gamma=0.99999$, $\delta=0.5682\cdot\frac{9}{D}$ and
$\varepsilon=10^{-9}$. Using $L_R\geq\frac{9}{D}q_C$ in the first
two inequalities in \eqref{eq:main-uniform-simple-certificates}
gives
\[
 S\leq\frac54-\frac{\gamma}{8}+2\varepsilon+\gamma q_C,
 \qquad
 S\leq\frac54+\varepsilon-\delta q_C.
\]
Multiply the first bound by $\delta$ and the second by $\gamma$,
then add and divide by $\gamma+\delta$. The terms in $q_C$
cancel, giving
\[
 S\leq\frac54-\frac{\gamma\delta}{8(\gamma+\delta)}
         +\frac{(2\delta+\gamma)\varepsilon}{\gamma+\delta}.
\]
This is less than $\frac{1211}{1000}<\frac{147}{121}=M$ when
$D=11$ and less than $\frac{1215}{1000}<\frac{206}{169}=M$
when $D=13$. Thus $b>0$ in every case, as also recorded for
$D\geq9$ in \textbf{Lemma~S4.3}
of~\cite{FGSTSupplement}.

To obtain degree bounds at a pair joined by parallel marked edges,
fix vertices $u,v$ joined by exactly $k\geq2$ marked edges. Call
$u,v$ the \emph{roots}. Let
\[
 A=(N_G(u)\cup N_G(v))\setminus\{u,v\},\qquad
 q=|A|,\qquad W=V(G)\setminus(A\cup\{u,v\}).
\]
Let $H$ be the spanning subgraph of $G$ consisting of its marked
edges. We call each edge between $\{u,v\}$ and $A$ a \emph{spoke}.
At most $D-k$ spokes are incident with each root, so
$q\leq2(D-k)$. Every marked edge other than the $k$ edges with
endpoints $u,v$ has an endpoint in $A$. An edge incident with
exactly one root is a spoke. An edge disjoint from the roots
has an endpoint in $A$ by the strong clique property.

The sum $d_H(u)+d_H(v)$ counts each of the $k$ joining edges
twice, so at most $2D-k$ marked edges meet the roots. Each
vertex of $A$ has at least one edge to a root and can therefore
be incident with at most $D-1$ edges avoiding the roots.
Summing over $A$ bounds the marked edges disjoint from the
roots by $(D-1)q$, giving the first line below. To count all the
marked edges that have an endpoint in $A$, sum $d_H(a)$ over
$a\in A$ and subtract one for each marked edge with both
endpoints in $A$. Adding the $k$ edges between the roots gives
the second line.
\begin{equation}\label{eq:odd-degree-root-count}
\begin{aligned}
 m&\leq 2D-k+(D-1)q\leq 2D^2-(2D-1)k,\\
 m&=k+\sum_{a\in A}d_H(a)-|E(H[A])|.
\end{aligned}
\end{equation}
These counts apply to any parallel marked pair and will be used
at $D=3,5,7,9$.

\subsection{Weighted Density Estimates for \texorpdfstring{$D\geq9$}{D at least 9}}
\label{subsec:exact-large-degree}

The first column of \textbf{Table~\ref{tab:multi-odd}} partitions
the odd integers $D\geq9$. The expression for $M$ in
\eqref{eq:odd-degree-mass}, as a function of real $D$, has
derivative $\frac{D-5}{2D^3}>0$ for $D\geq9$. For each finite
degree range, evaluate $M$ at its two endpoints and take a
rational interval $[M_-,M_+]$ containing both values. For the
range $D\geq59$, take $M_+=\frac54$, since
\eqref{eq:odd-degree-mass} gives $M<\frac54$ on this range.

Let $\alpha$ be a lower bound on $\frac{R(C)}{D}$ that holds
for every induced red five-cycle at every degree in a given row.
\textbf{Lemma~S4.6} of~\cite{FGSTSupplement} gives
$R(C)\geq11$ and $2\leq b\leq17$ at $D=9$.
\textbf{Lemma~S4.2} of~\cite{FGSTSupplement}
gives $R(C)\geq9$ at $D=11$ and $R(C)\geq11$ at $D=13,15$.
For the other odd degrees, $R(C)>\frac{D+1}{2}$ and the parity
of $R(C)$ give
\[
 R(C)\geq
 \begin{cases}
  \frac{D+5}{2},&D\equiv1\pmod4,\\
  \frac{D+3}{2},&D\equiv3\pmod4 .
 \end{cases}
\]
For example, at $D=21,23$ the bound is $R(C)\geq13$, so
$\frac{R(C)}{D}\geq\frac{13}{23}>0.565$. The table records a
valid $\alpha$ for each row.

Summing $R(C)\geq\alpha D$ over the ordered cycles gives
$L_R\geq\alpha q_C$. This supplies the weighted density constraint
\eqref{eq:main-uniform-cycle-constraint}. Substituting the same
bound in \eqref{eq:main-uniform-simple-certificates} and using
$S\leq\frac54$ gives $S\leq S_{\max}(\alpha)$ by the calculation
in \textbf{Subsection~\ref{subsec:uniform-weighted-constraints}}.
Since $e=1-\frac{S}{M}$, it follows that
\[
 e\geq 1-\frac{S_{\max}(\alpha)}{M_-}.
\]
The rows with $D\geq15$ use this lower bound on $e$. At $D=9$,
the bounds on $b$ give
$\frac{2}{98}\leq e\leq\frac{17}{98}$. At $D=11,13$,
\textbf{Lemma~S4.5} of~\cite{FGSTSupplement}
excludes $b=1$. Since $b>0$, this gives $e\geq\frac{2}{m}$.

To bound $e$ above on the other finite ranges, let $j_v$ count
the remaining marked edges incident with $v$ and put
$a=\frac{\sum_v j_v^2}{2bD}$. \textbf{Lemma~S3.12}
of~\cite{FGSTSupplement} gives
\[
 M\leq2-\left(2+\frac{1}{D}\right)a-\frac{1}{D^2},
 \qquad Me\leq2a(1-a)-\frac{a}{D},
 \qquad e\leq1-\frac{M}{2}.
\]
If $D_+$ is the largest degree in the row, then $D\leq D_+$
and $M\geq M_-$ give
\[
 a\leq a_+:=\frac{2-M_--D_+^{-2}}{2+D_+^{-1}}.
\]
In these rows $M_->\frac65$, so $a_+\leq\frac{2-M_-}{2}<\frac25$.
Since $D_+\geq9$, we have $\frac25<\frac12-\frac{1}{4D_+}$.
Thus $2a(1-a)-\frac{a}{D_+}$ is increasing on $[0,a_+]$.
Also $\frac{a}{D}\geq\frac{a}{D_+}$. The second inequality
therefore gives
\[
 e\leq\frac{2a_+(1-a_+)-\frac{a_+}{D_+}}{M_-}.
\]
We use this bound or $e\leq1-\frac{M_-}{2}$, whichever is
smaller. For $D\geq59$, \eqref{eq:odd-degree-mass} gives
$M\geq\frac{4323}{3481}>\frac{4967}{4000}$.
\textbf{Lemma~S3.11} of~\cite{FGSTSupplement} then gives $e<\frac{1}{4}$,
which supplies the upper endpoint of the $e$-interval in the last row
of \textbf{Table~\ref{tab:multi-odd}}.

\begin{table}[t]
\centering
\small
\begin{tabular}{@{}lcccc@{}}
\toprule
Odd $D$ & $[M_-,M_+]$ & $[e_-,e_+]$ & $\alpha$ & Verified conclusion\\
\midrule
$9$        & $[1.209,1.210]$  & $[0.020,0.180]$   & $1.222$ & $M<-1.973$\\
$11$       & $[1.2148,1.215]$ & $[0.0135,0.360]$  & $0.818$ & $M<0.785$\\
$13$       & $[1.2185,1.219]$ & $[0.009,0.362]$   & $0.846$ & $M<0.573$\\
$15$       & $[1.222,1.223]$  & $[0.007,0.365]$   & $0.733$ & $M<0.827$\\
$17$       & $[1.2249,1.225]$ & $[0.006,0.368]$   & $0.647$ & $M<0.631$\\
$19$       & $[1.227,1.228]$  & $[0.006,0.370]$   & $0.578$ & $M<0.910$\\
$21$--$23$ & $[1.2249,1.231]$ & $[0.004,0.372]$   & $0.565$ & $M<1.227683$\\
$25$--$37$ & $[1.232,1.238]$  & $[0.009,0.384]$   & $0.542$ & $M<0.177$\\
$39$--$47$ & $[1.238,1.240]$  & $[0.0134,0.382]$  & $0.531$ & $M<-0.676$\\
$49$--$57$ & $[1.240,1.242]$  & $[0.0082,0.380]$  & $\frac{1}{3}$ & $M<0.292$\\
$\ge59$    & $[1.24175,1.25]$ & $[0.005,0.250]$   & $0.228$ & infeasible\\
\bottomrule
\end{tabular}
\caption{Parameter intervals and verified bounds for odd $D$.
All decimal entries denote rational numbers.}
\label{tab:multi-odd}
\end{table}

To exclude each row, apply the extension equations
\eqref{eq:weighted-extension}, the weighted flag inequalities,
and the linear constraints from the bounds on $M$, $e$ and $R(C)$.
Every pattern containing a red $K_5$ has density zero. Indeed,
at $m=f(D)+1$, we have $4m=5D^2-2D+5$, which exceeds the bound
for a red $K_5$ in (S3.27) of~\cite{FGSTSupplement}.
The rows with $D\leq23$ also use the absence of four vertices
whose six pairs are red, by \textbf{Lemma~S4.7}
of~\cite{FGSTSupplement}. The proof of \textbf{Lemma~S4.6}
of~\cite{FGSTSupplement} gives $13\leq q\leq14$ in the reduced
multigraph for every parallel marked pair at $D=9$. The regular
completion preserves those neighbors and the marked edges.
Since each root has two marked edges to the other root, each
has at most seven distinct neighbors outside the root pair
in the completion. Thus the union of their neighborhoods outside
the root pair still has size between $13$ and $14$. To count the
possible images of an added position adjacent to at least one
root image, take the union of the full neighborhoods. That union
also contains the two roots, giving between $15$ and $16$ images.
At $D=13,15$, the
certificates use the inequalities
obtained from the possible red triangle neighborhoods in
\textbf{Figure~S4.17} of~\cite{FGSTSupplement}.
\textbf{Appendix~SG.3} of~\cite{FGSTSupplement}
gives the enumeration and explains how these restrictions enter
the linear constraints.

\textbf{Appendix~SG.4} of~\cite{FGSTSupplement}
verifies the multipliers and the correction for negative residual
coefficients for the rows of \textbf{Table~S4.1} there.
On every finite degree range the
rounded upper bound for $M$ in the last column is below the
smallest possible $M$ from \eqref{eq:odd-degree-mass}. For
$D=21,23$, for example, the smallest possible value is
$\frac{542}{441}>1.227683$.

For $D\geq59$, take the linear form to be bounded to be zero
in place of $M$ in \eqref{eq:main-uniform-final-decomposition}.
For any density vector satisfying the constraints in this row,
the resulting identity expresses $C+B$ as a sum of three
nonnegative terms. Here $C$ is the multiplier of the normalization
equation and $B$ bounds $\sum_{j:r_j<0}(-r_j)x_j$, using the
residual coefficients $r_j$ and coordinates $x_j$ for this row.
Equation (SG.3) of~\cite{FGSTSupplement} gives $C+B<0$,
so the row has no feasible density vector. This excludes
every row in \textbf{Table~\ref{tab:multi-odd}}.

\subsection{The Cases \texorpdfstring{$D=1,3,5,7$}{D = 1, 3, 5, 7}}
\label{subsec:exact-small-degree}

At $D=1$ every component has at most one edge, so any strong
clique has size at most $1=f(1)$. For $D=3,5,7$ use the roots
$u,v$ and the sets $A,W$ of
\textbf{Subsection~\ref{subsec:exact-structure}}. When $D=3$,
\eqref{eq:odd-degree-root-count} gives $m\leq18-5k\leq8$,
contrary to $m=f(3)+1=11$.

Suppose $D=5$, so $m=f(5)+1=30$. The first line of
\eqref{eq:odd-degree-root-count} gives $30\leq50-9k$ and
therefore $k=2$. It then gives $30\leq8+4q$, while
$q\leq2(5-k)=6$, so $q=6$. Since each root has at most
three spokes and $|A|=6$, each root has three distinct neighbors
in $A$, the neighbor sets are disjoint, and the spokes have
multiplicity one. Let $p$ count the unmarked spokes and let
$i_m,i_u$ count the marked and unmarked edges within $A$.
The $p$ unmarked spokes contribute one unmarked incidence each
at a vertex in $A$, and the $i_u$ unmarked edges within $A$
contribute two each.
The second line of \eqref{eq:odd-degree-root-count} gives
$\sum_{a\in A}d_H(a)=28+i_m$. Since each of the six vertices
of $A$ has degree at most five,
\[
 p+i_m+2i_u \leq
 \sum_{a\in A}\bigl(5-d_H(a)\bigr)+i_m=2.
\]
There are 21 arrangements of the marked and unmarked edges on
$A\cup\{u,v\}$ satisfying this inequality, up to permuting the
neighbors of either root and interchanging the roots.
\textbf{Figure~S4.21} of~\cite{FGSTSupplement}
lists them.

Fix an arrangement on $A\cup\{u,v\}$ and write
$A=\{a_1,\ldots,a_6\}$. Let $\tau_i\in\{u,v\}$ be the root
adjacent to $a_i$. The number of edges from $a_i$ to $W$ is at
most $c_i=4-d_{G[A]}(a_i)$, since one edge joins $a_i$ to
$\tau_i$. To count the marked edges from a vertex $x\in W$
to $A$, let $z_i$ be the number of marked edges with endpoints
$x,a_i$. Write $z=(z_i)_{i=1}^6$. For each $j$, let $I_j(z)$
consist of the indices $i$ for which $z_i>0$, $\tau_i\ne\tau_j$,
$a_j\tau_j$ is marked, and $a_ia_j\notin E(G)$.
If $i\in I_j(z)$, then $xa_j\in E(G)$. Indeed, a
marked edge with endpoints $x,a_i$ and the marked spoke
$a_j\tau_j$ require a joining edge. Of the four possible
joining pairs, $x\tau_j$, $a_i\tau_j$ and $a_ia_j$ are
nonadjacent, leaving $xa_j$.

Put
\[
 r_j=
 \begin{cases}
 z_j,&I_j(z)=\varnothing,\\
 \max\{z_j,1\},&I_j(z)\ne\varnothing.
 \end{cases}
\]
The number of edges with endpoints $x,a_j$ is
at least $r_j$. Thus $\sum_jr_j\leq5$ by the degree bound at
$x$. Also $\sum_i z_i>0$ because every vertex of the reduced
multigraph is incident with a marked edge and every marked edge
with an endpoint in $W$ has its other endpoint in $A$.

To bound the total number of marked edges between $A$ and $W$,
write $c=(c_i)_{i=1}^6$. For each $c'\in\mathbb Z_{\geq0}^6$
with $c'\leq c$, use $c'_i$ as an upper bound on the number of
edges from $a_i$ to vertices of $W$ still to be considered.
For each $z=(z_i)_{i=1}^6\in\mathbb Z_{\geq0}^6$ with
$z\leq c'$, compute $r=(r_j)_{j=1}^6$ by the preceding rule.
Define
\[
 V(c')=\max\left(
  \{0\}\cup
  \left\{\sum_i z_i+V(c'-r):
       \sum_i z_i>0,\ \sum_i r_i\leq5,\ r\leq c'\right\}
 \right),
\]
where vector inequalities are coordinatewise. Each recursive call
decreases the sum of the coordinates of $c'$, since
$\sum_i r_i\geq\sum_i z_i>0$. Consider the vertices of $W$ in any
order. For each vertex, its vector $z$ is among the choices in the
maximum because the number of edges used at $a_i$ is at least
$r_i$. Starting from $c$ and subtracting $r$ at each step therefore
leaves an upper bound on the edges at each $a_i$ to the remaining
vertices. Thus $V(c)$ bounds the total number of marked edges
between $A$ and $W$. The degree-five computation in
\textbf{Appendix~SF.1} of~\cite{FGSTSupplement} verifies the recurrence
for each arrangement.
Adding the marked edges fixed on $A\cup\{u,v\}$ to $V(c)$
gives at most $27$ for each arrangement, contrary to $m=30$.

Suppose $D=7$, so $m=f(7)+1=59$. The inequality
$59\leq98-13k$ excludes $k\geq4$. If $k=3$, then $q\leq8$ and
the second line of \eqref{eq:odd-degree-root-count} gives
\[
 59=3+\sum_{a\in A}d_H(a)-|E(H[A])|\leq3+7q\leq59.
\]
Equality forces $q=8$, $d_H(a)=7$ for all $a\in A$, and
$E(H[A])=\varnothing$. Since $d_G(a)\leq7$, every edge
incident with $A$ is marked and $G[A]$ has no edges. For
$z\in\{u,v\}$, let $A_z=N_G(z)\setminus\{u,v\}$. Each root
has at most four spokes, so $A_u,A_v$ are disjoint, each has
size four, and each vertex in $A$ has one spoke. Fix
$a\in A_u$. Because $d_H(a)=7$ and $G[A]$ has no edges, there
is a marked edge $ax$ with $x\in W$. For every $b\in A_v$,
the marked edges $ax$ and $vb$ require a joining edge. The
pairs $av$, $ab$ and $xv$ are nonadjacent, so $xb$ is an edge.
It is marked because $d_H(b)=7$. Fix $b\in A_v$. For each
$a'\in A_u$, the marked edges $xb$ and $ua'$ similarly
require the adjacency $xa'$. Thus $x$ is adjacent to all
eight vertices of $A$, a contradiction.
Thus every pair of vertices joined by parallel marked edges
has marked multiplicity two.

Now fix a pair of marked multiplicity two. The inequality
$59\leq12+6q$ and the bound $q\leq10$ give $q\in\{9,10\}$.
Let $T$ count the spokes with multiplicity and $p$ count the
unmarked spokes. Write $i_m,i_u$ for the marked and unmarked
edges within $A$, and $m_\times,U$ for the marked and unmarked
edges between $A$ and $W$, respectively. Put
$\ell=\sum_{a\in A}(7-d_G(a))$. The marked-edge count and the
degree sum over $A$ give
\[
 57=T-p+i_m+m_\times,\qquad
 \ell+p+i_m+2i_u+U=7q-57.
\]
The first equality counts the marked edges other than the two
between the roots. For the second, write
$\sum_{a\in A}d_G(a)=T+2(i_m+i_u)+m_\times+U=7q-\ell$
and substitute for $m_\times$ from the first equality.
The second equality gives $i_m\leq7q-57-p$. Substituting
this upper bound for $i_m$ in the first equality and using
$T,q\leq10$ gives $m_\times\geq114-T-7q+2p\geq34+2p$.

To bound $|W|$, put $w=|W|$ and let $D_A=m_\times+U$ be the
number of edges between $A$ and $W$. For each $x\in W$, choose
$a(x)\in A$ joined to $x$ by a marked edge. This is possible
because every vertex of $G$ is incident with a marked edge
and all marked edges from $W$ have an endpoint in $A$. Write
$n_a=|\{x\in W:a(x)=a\}|$ and let $G_A$ be the underlying
simple graph of $G[A]$. We count unordered pairs $\{x,y\}$
for which $a=a(x)$ and $b=a(y)$ are distinct and
$ab\notin E(G_A)$. For each such pair, the marked edges $ax,by$
require at least one of $xy,ay,bx$ as a joining edge. Hence
\[
 \frac{w^2-\sum_a n_a^2}{2}
     -\sum_{ab\in E(G_A)}n_an_b
 \leq |E(G[W])|+\sum_a n_a\bigl(d_G(a,W)-n_a\bigr),
\]
The left-hand side counts all pairs in different sets
$\{x:a(x)=a\}$ and subtracts those whose chosen neighbors are
adjacent in $G_A$. An edge in $G[W]$ accounts for at most one
pair. An edge from $a\in A$ to $y\in W$ with $a(y)\ne a$ accounts
for at most $n_a$ pairs, one for each $x$ with $a(x)=a$. The
number of these edges at $a$ is at most $d_G(a,W)-n_a$, where
$d_G(a,W)$ counts edges with multiplicity. Each $a$ has at
least one spoke and at least $d_{G_A}(a)$ edges within $A$, so
$d_G(a,W)\leq6-d_{G_A}(a)$. The degree sum over $W$ gives
$2|E(G[W])|\leq7w-D_A$. Substituting the two bounds and
rearranging gives
\[
  w^2+\sum_a n_a^2+2\sum_{ab\in E(G_A)}(n_a+n_b-n_an_b)
 \leq19w-D_A.
\]
Using $2n_an_b\leq n_a^2+n_b^2$ for each edge $ab\in E(G_A)$,
the left side is at least
$w^2+\sum_a[n_a^2+d_{G_A}(a)n_a(2-n_a)]$.
For $n_a\leq2$, each term in the brackets is nonnegative.
For $n_a\geq3$, use $d_{G_A}(a)\leq6-n_a$ and $2-n_a<0$ to
bound the bracket from below by $n_a(n_a-3)(n_a-4)$, which
is nonnegative for integer $n_a\geq3$. Since $D_A\geq34$,
we have $w^2\leq19w-34$ and hence $w\leq17$.
If $w=17$, equality forces all positive $n_a$ to be three
or four, with $d_{G_A}(a)=3$ when $n_a=3$. Equality in
$2n_an_b\leq n_a^2+n_b^2$ also forces $n_a=n_b$ on every edge
$ab\in E(G_A)$. The only multiset of threes and fours summing
to $17$ is $\{3,3,3,4,4\}$. Thus each vertex with $n_a=3$ can
have at most two neighbors in $G_A$, contradicting
$d_{G_A}(a)=3$.
It follows that $w\leq16$.

To bound $m_\times$ from above, take the induced host on $A\cup W$
and consider only the marked edges between $A$ and $W$.
These edges form a bipartite subgraph and remain a strong clique
in the induced host, since their endpoints are in $A\cup W$.
The host degrees on
$A$ and $W$ are at most six and seven, respectively, so the marked
subgraph has Ore-degree at most 13 in this host. The finite recurrence
based on (S4.14) of~\cite{FGSTSupplement} adapts the neighborhood count
of Kumar, Mohar and Pragada~\cite{KumarMoharPragada2026} to parallel
marked edges. For $D=7$, the computation in \textbf{Appendix~SF.1}
of~\cite{FGSTSupplement}
gives $m_\times\leq48$ at Ore-degree at most 13 and
$m_\times\leq41$ at Ore-degree at most 12. In particular,
$114-T-7q+2p\leq48$ gives $p\leq7$ when $q=10$
and $p\leq3$ when $q=9$.

Call a neighbor of a root in $A$ \emph{private} if it is
adjacent to exactly one root. For $q=10$, each root has five
private neighbors and all spokes have multiplicity one. Order
the roots so that their marked spoke counts $a,b$ satisfy $a\leq b$.
Since $p\leq7$, we have $a+b=10-p\geq3$. Of the
21 pairs with $0\leq a\leq b\leq5$, the four pairs
$(0,0),(0,1),(0,2),(1,1)$ have sum less than three, leaving
17 choices.
For $q=9$ with disjoint root neighborhoods in $A$, one root
has four neighbors in $A$ and the other has five, since each
root has at most five spokes. The first root has either four
spokes of multiplicity one or two parallel marked spokes to
one neighbor and one spoke to each of its three other neighbors.
The second root has five spokes of multiplicity one. In either case
the numbers of unmarked spokes at the two roots can be any
ordered pair of nonnegative integers with sum at most three.
This gives $1+2+3+4=10$ choices in each case.
For $q=9$ with a common neighbor in $A$, both roots must have
five neighbors in $A$, of which exactly one is common. Every
spoke has multiplicity one. If both spokes
to the common neighbor are unmarked, the private unmarked
spoke counts sum to at most one, giving two choices up to
root interchange. If exactly one of the common spokes is
unmarked, fix its root first. The two private counts can then
be any ordered pair with sum at most two, giving six choices.
If both common spokes are marked, the private counts sum to
at most three. There are six choices up to root interchange.

To exclude the $17+10+10+14=51$ arrangements, fix the two
marked edges between the roots and all spokes for each
arrangement. Assign 16 vertices to $W$ in the Boolean formula.
A graph with fewer than 16 vertices in $W$ can be represented
by setting all the variables at the unused vertices to zero.
For each unordered pair
$\{x,y\}$, the Boolean variables $E_{xy}$, $M_{xy}$ and $J_{xy}$
record whether $x,y$ are adjacent, whether at least one marked
edge joins them, and whether two marked edges join them,
respectively. The formula enforces
$J_{xy}\Rightarrow M_{xy}\Rightarrow E_{xy}$. In the reduced
host, an adjacent pair with no marked edge has exactly one
unmarked edge. A pair with a marked edge has only marked edges,
and its multiplicity is at most two. Hence the host and marked
multiplicities on the pair are $E_{xy}+J_{xy}$ and
$M_{xy}+J_{xy}$, respectively, with Boolean values read as $0$ or $1$.
The formula imposes the degree bounds, a total of $59$ marked
edges, the joining clause
(SF.1) of~\cite{FGSTSupplement} for
each disjoint pair of candidate marked edges, and the
necessary capacity and marked-edge count bounds just derived.
It also imposes the bounds from the Ore-degree recurrence and
\eqref{eq:odd-degree-root-count} at each pair of vertices
joined by two marked edges. A counterexample satisfies all
these conditions, so its host and marked multiplicities would
give a satisfying assignment for its arrangement. The 51
arrangements are given in
\textbf{Figures~S4.22--S4.24} of~\cite{FGSTSupplement}.
\textbf{Appendix~SF.1} of~\cite{FGSTSupplement} verifies their
refutations by reverse unit propagation. In this check,
an added clause is justified by deriving a contradiction from its
negation by unit propagation.
The check verifies $22{,}217{,}393$ clause additions in total.
Every refutation derives the empty clause, excluding
the corresponding arrangement. This completes the proof of
\textbf{Theorem~\ref{thm:multigraph}}.

\section*{Acknowledgement}
This research was supported in part by the National Natural Science Foundation
of China (Nos. 12471330 and 12571373) and the Youth Student Fundamental Research
Funds of Shandong University.

\section*{Declaration on the Use of Generative AI}
The authors developed several key ideas in the proof, including
the weighted flag algebra method and the strategy of analyzing
the structure of the hypothetical counterexample of
\textbf{Conjecture~\ref{conj:fgst}}, while we used ChatGPT to assist
with the exploration of some proof strategies and calculations
and with auxiliary tasks, including improving the organization
and clarity of the manuscript and addressing some possible
inconsistencies or gaps.
\relax\space

\begingroup
\small
\setlength{\bibsep}{3pt plus .3pt minus .3pt}
\phantomsection
\addcontentsline{toc}{section}{\refname}
\bibliography{refs}
\endgroup
\end{document}